\documentclass{elsart}
\usepackage{amsmath}
\usepackage{amsfonts}
\usepackage{amssymb}
\usepackage{subfigure}
\usepackage{graphicx}
\usepackage{natbib}

\theoremstyle{definition}
\newtheorem{assumption}{Assumption}

\newcommand{\psca}[1]{\langle #1 \rangle}	
\newcommand{\abs}[1]{\lvert#1\rvert}		
\newcommand{\norm}[1]{\lVert#1\rVert}		
\newcommand{\normsup}[1]{\norm{#1}_{\infty}}	

\newcommand{\blackboard}[1]{\mathbb{#1}} 
\renewcommand{\P}{\blackboard{P}}		
\newcommand{\E}{\blackboard{E}}			
\renewcommand{\L}{\blackboard{L}}		

\newcommand{\calig}[1]{\mathcal{#1}}
\newcommand{\A}{\calig{A}} 		
\newcommand{\C}{\calig{C}} 		
\newcommand{\F}{\calig{F}} 		

\DeclareMathOperator{\Tr}{Tr}		
\DeclareMathOperator{\Id}{Id}		
\DeclareMathOperator{\var}{var}

\begin{document}

\begin{frontmatter}
\title{An adaptive scheme for the approximation of dissipative systems}
\author{Vincent Lemaire}
\ead{vincent.lemaire@univ-mlv.fr}
\address{Laboratoire d'Analyse et de Math\'{e}matiques Appliqu\'{e}es, UMR 8050,\\ Universit\'{e} de Marne-la-Vall\'{e}e, 5 boulevard Descartes, Champs-sur-Marne, \\ F-77454 Marne-la-Vall\'{e}e Cedex 2, France.} 
\begin{abstract}
We propose a new scheme for the long time approximation of a diffusion when the drift vector field is not globally Lipschitz. Under this assumption, regular explicit Euler scheme --with constant or decreasing step-- may explode and implicit Euler scheme are CPU-time expensive. The algorithm we introduce is explicit and we prove that any weak limit of the weighted empirical measures of this scheme is a stationary distribution of the stochastic differential equation. Several examples are presented including gradient dissipative systems and Hamiltonian dissipative systems. 
\end{abstract}
\end{frontmatter}

\section{Introduction}
We consider the following stochastic differential equation  
\begin{equation} \label{intro-diffusion}
	dx_t = b(x_t) dt + \sigma(x_t) dB_t, \qquad x(0) = x_0 \in \Rset^d,
\end{equation}
where $b:\Rset^d \rightarrow \Rset^d$ is a locally Lipschitz continuous vector field and $\sigma$ is locally Lipschitz continuous on $\Rset^d$, with values in the set of $d \times m$ matrices and $B$ is an $m$-dimensional Brownian motion. Assume that $(x_t)_{t \ge 0}$ has a Lyapounov function $V$ $i.e.$ a positive regular function decreasing along trajectories (precise conditions are given by Assumption~1 in Section~\ref{section:framework}), so that there exists at least one invariant measure.

Until recently, the approximation of the stationary mode of the diffusion has been studied under the assumption that $V$ is \emph{essentially quadratic} $i.e.$ 
\begin{equation} \label{hypotheses-LambPag}
	\abs{\nabla V}^2 = \calig{O}(V) \quad \text{and} \quad \sup_{x \in \Rset^d} \norm{D^2 V} < + \infty, 
\end{equation}
and $\abs{b}^2 = \calig{O}(V)$ (which implies sublinear growth for $b$). 
When $\sigma$ is bounded and the diffusion is uniformly strictly elliptic, the invariant measure $\nu$ is unique and Talay proposed in \cite{talay-90} a method for the computation of $\nu$ based on the \emph{constant step} Euler scheme. He proved the convergence of the invariant measure of the scheme to $\nu$. On the other hand, Lamberton and Pag\`{e}s studied in \cite{lamberton-pages-02} the ergodic properties of the weighted empirical measures $(\nu^\eta_n)_{n \ge 0}$ of a \emph{decreasing step} Euler scheme. They proved the almost sure tightness of $(\nu_n^{\eta})_{n \ge 1}$ and that any weak limit is a stationary distribution for the diffusion.

However, the conditions \eqref{hypotheses-LambPag} and $\abs{b}^2 = \calig{O}(V)$ are too restrictive for studying systems used in random mechanics (see Soize \cite{soize}). Indeed, the drift vector field $b$ is generally locally Lipschitz and in many cases $V$ is not essentially quadratic. This framework has been recently investigated by Talay in \cite{talay-02} and by Mattingly et al. in \cite{mattingly-02}. In these papers, \emph{implicit} Euler schemes with constant steps are used for the approximation of the diffusion. In recent work, Lamba, Mattingly and Stuart have introduced on finite time interval $[0; T]$ an adaptive explicit Euler scheme (see \cite{stuart-ad-part1} and \cite{stuart-ad-part2}). The step is adapted according to the error between the Euler and Heun approximations of the ODE $\dot{x} = b(x)$. They prove strong mean-quadratic convergence of the scheme on over finite time intervals and ergodicity when the noise is non-degenerate. We propose a completely different \emph{explicit scheme based on a stochastic step sequence} and we obtain the almost sure convergence of its weighted empirical measures to the invariant measure of \eqref{intro-diffusion}.

The key to prove the almost sure tightness of the weighted empirical measures of the decreasing Euler scheme $(X_n)_{n \ge 0}$ introduced in \cite{lamberton-pages-02} is that the scheme satisfies a \emph{stability condition} $i.e.$ there exist $\tilde{\alpha} > 0$ and $\tilde{\beta} > 0$ such that 
\begin{equation} \label{condition-rappel-algo}
	\frac{\E \bigl[ V(X_{n+1}) \,|\, \F_n \bigr] - V(X_n)}{\gamma_{n+1}} \le - \tilde{\alpha} V(X_n) + \tilde{\beta},
\end{equation}
where $(\gamma_n)_{n \ge 0}$ is the deterministic decreasing step sequence. Without assumptions \eqref{hypotheses-LambPag} we can no longer prove the stability condition \eqref{condition-rappel-algo} for this scheme. Our scheme is built in order to satisfy \eqref{condition-rappel-algo}. We proceed as follows. Firstly, we start from a deterministic $X_0 = x_0 \in \Rset^d$ and set 
\begin{equation} \label{algorithme}
	X_{n+1} = X_n + \tilde{\gamma}_{n+1} b(X_n) + \sqrt{\tilde{\gamma}_{n+1}} \sigma(X_n) U_{n+1}, \quad n \ge 0,
\end{equation}
where $(U_n)_{n \ge 1}$ is a $\Rset^m$-white noise more precisely defined in Section~\ref{section:framework} and $\tilde{\gamma}_{n+1} = \gamma_{n+1} \wedge \chi_n$ with $(\gamma_n)_{n \ge 0}$ a positive nonincreasing sequence and $\chi_n$ a $\sigma(U_1, \dots, U_n)$--measurable random variable.   

The basic main idea is to choose $\chi_n$ small when the scheme starts to explode. In this case the discretization is finer and the stability condition of the diffusion prevents the explosion. Furthermore we prove that the scheme satisfies a similar condition with \eqref{condition-rappel-algo}. A non-optimal --although natural-- choice for $\chi_n$ may be  $\chi_n = \frac{1}{\abs{b(X_{n-1})}^2 \vee 1}$. For the two studied examples, optimal choices depend on the Lyapounov function (see \eqref{thm-essquad-majo-chi} and \eqref{thm-hamilton-majo-chi}).

A crucial feature of our algorithm is the existence of an almost surely finite time $n_1$ such that for every $n \ge n_1$, $\tilde{\gamma}_n = \gamma_n$ i.e. the event $\{\chi_{n-1} < \gamma_n\}$ does not occur any more. 

Numerically the algorithm is very simple to implement and the complexity is the same as that of a regular Euler scheme. Another interest is that the scheme is explicit, which is a big advantage on implicit schemes for high dimensional problems. Indeed a fixed point algorithm is not needed is our case. Moreover, we will see that wrong convergence problem due to fixed point algorithm may be avoided using our algorithm. 

The paper is organized as follows. We introduce the framework and the algorithm in Section~2. In Section~3 are presented some preliminary results about the approximation scheme of $(X_n)_{n \ge 0}$ defined in \eqref{algorithme}. In Section~4 we extend some results of \cite{lamberton-pages-02} and give conditions for the almost sure tightness of the empirical measure and for its weak convergence to an invariant measure of \eqref{intro-diffusion}. Section~5 is devoted to the study of monotone systems 
and Section~6 of stochastic Hamiltonian dissipative systems. The numerical experiments are in Section~7 including some comparaison with recently introduced implicit scheme. We confirm the non-explosion and the convergence of the scheme. 

\section{Framework and algorithm} \label{section:framework}
We will denote by $\A$ the infinitesimal generator of \eqref{intro-diffusion}. 
The following assumption will be needed throughout the paper.
\begin{assumption} 
There is a $\C^2$ function $V$ on $\Rset^d$ with values in $[1, +\infty[$ such that $\displaystyle \lim_{\abs{x} \rightarrow +\infty} V(x) = +\infty$ and satisfying 
	\begin{equation} \label{condition-rappel}
		\exists \alpha > 0, \; \exists \beta > 0, \quad \psca{ \nabla V, b } \le -\alpha V + \beta,
	\end{equation}
	\begin{gather} \label{condition-sigma-V}
	\begin{split}
		\exists C_{V, \sigma} > 0, \; \exists a \in (0, 1], \quad \Tr \bigl( \sigma^* (\nabla V)^{\otimes 2} \sigma \bigr) &\le C_{V, \sigma} V^{2-a}, \quad \\
		\text{and} \qquad \exists C > 0, \quad \sup_{x \in \Rset^d} \Tr \bigl( \sigma^* D^2V \sigma \bigr)(x) &\le C \sup_{x \in \Rset^d} \Tr (\sigma^* \sigma)(x).
		\end{split}
	\end{gather}
\end{assumption}
\begin{rem}
If $V$ is essentially quadratic ($i.e.$ satisfies \eqref{hypotheses-LambPag}) then \eqref{condition-sigma-V} is satisfied as soon as there exists $C_\sigma > 0$ and $a \in (0,1]$ such that $\Tr(\sigma^* \sigma) \le C_\sigma V^{1-a}$.
\end{rem}

Under this assumption, there exists a global solution to equation \eqref{intro-diffusion} and (at least) one invariant measure. An important point to note here is that all invariant measures have exponential moments.
Indeed, an easy computation shows that for all $\lambda < \frac{\alpha}{a C_{V, \sigma}}$ we have 
\begin{equation*}
	\exists \tilde{\alpha} > 0, \; \exists \tilde{\beta} > 0, \quad \A \exp \bigl( \lambda V^a \bigr) \le - \tilde{\alpha} \exp(\lambda V^a) + \tilde{\beta}, 
\end{equation*}
and this implies $\nu \bigl( \exp(\lambda V^a) \bigr)$ is finite for all invariant measures $\nu$. 

For the approximation of the diffusion, we assume that $(U_n)_{n \ge 1}$ is a sequence of i.i.d. random variables defined on a probability space $(\Omega, \mathfrak{A}, \P)$, with values in $\Rset^m$, and such that $U_1$ is a generalized Gaussian (see Stout \cite{stout}) $i.e.$ 
\begin{equation} \label{bruit-sous-gaussien}
	\exists \kappa > 0, \; \forall \theta \in \Rset^d, \quad \E \bigl[ \exp \bigl( \psca{\theta, U_1} \bigr) \bigr] \le \exp \Bigl( \frac{\kappa \abs{\theta}^2}{2} \Bigr).
\end{equation}
and that $\var(U_1) = \Id_m$. We will call $(U_n)_{n \ge 1}$ a $\Rset^m$-valued generalized Gaussian white noise. The condition \eqref{bruit-sous-gaussien} implies that $U_1$ is centered and satisfies 
\begin{equation} \label{bruit-moment-exponentiel}
	\exists \tau > 0, \quad \E \bigl[ \exp \bigl( \tau \abs{U_1}^2 \bigr) \bigr] < + \infty.
\end{equation}
Moreover, the condition $\var(U_1) = \Id_m$ implies that $\kappa \ge 1$.
In the sequel, $\F_n$ denotes, for $n \ge 1$, the $\sigma$-field generated on $\Omega$ by the random variables $U_1, \dots, U_n$, and $\F_0$ the trivial $\sigma$-field. 
\begin{rem}
	The assumptions made on the white noise $(U_n)_{n \ge 1}$ are not restrictive for numerical implementation. Indeed, centered Gaussian and centered bounded random variables satisfy \eqref{bruit-sous-gaussien}. 
\end{rem}

The \emph{stochastic step} sequence $\tilde{\gamma} = (\tilde{\gamma}_n)_{n \ge 0}$ is defined by 
\begin{equation} \label{definition-gamma-tilde}
	\forall n \ge 1, \quad \tilde{\gamma}_n = \gamma_n \wedge \chi_{n-1}, \qquad \tilde{\gamma}_0 = \gamma_0,
\end{equation}
where $(\gamma_n)_{n \ge 0}$ is a \emph{deterministic nonincreasing} sequence of positive numbers satisfying 
\begin{equation*}
	\lim_n \gamma_n = 0 \quad \text{and} \quad \sum_{n \ge 0} \gamma_n = + \infty,
\end{equation*}
and $(\chi_n)_{n \ge 1}$ is an $(\F_n)_{n \ge 0}$--adapted sequence of positive \emph{random} variables. It is important to note that the step sequence $\tilde{\gamma}$ is $(\F_n)_{n \ge 0}$--\emph{predictable}.

Now we introduce the \emph{weighted empirical measures} like Lamberton and Pag\`{e}s in \cite{lamberton-pages-02}. Given a sequence $\eta = (\eta_n)_{n \ge 1}$ of positive numbers satisfying $\sum_{n \ge 1} \eta_n = + \infty$, we denote by $\nu^{\eta}_n$ the random probability measure on $\Rset^d$ defined by 
\begin{equation*}
	\nu^{\eta}_n = \frac{1}{H_n} \sum_{k = 1}^n \eta_k \delta_{X_{k-1}}, \quad \text{with} \quad H_n = \sum_{k=1}^n \eta_k.
\end{equation*}
Throughout the paper, $\abs{.}$ denotes the Euclidean norm and $\norm{.}$ denotes the natural matrix norm induced by $\abs{.}$ $i.e.$ for every square matrix $A$, $\norm{A} = \sup_{\abs{x} = 1} \abs{A x}$. The letter $C$ is used to denote a positive constant, which may vary from line to line.

\section{Preliminary results}
In this section, we prove results which are the keys to study the Euler scheme with predictable random step defined in the introduction. Proposition \ref{lem-clef} contains two results: the first one \eqref{lem-clef-controle-esp} provides a substitute for the $L^p$--bounded-ness of $(V(X_n))_{n \ge 0}$ used in \cite{lamberton-pages-02}. The second one \eqref{lem-clef-controle-vitesse} is a new consequence of the stability condition \eqref{lem-condition-rappel-algo} and is used to prove the fundamental proposition \ref{prop-fondamentale} which ensures the existence of an almost surely finite time $n_1$ such that for every $n \ge n_1$, $\tilde{\gamma}_n = \gamma_n$. 

\begin{prop} \label{lem-clef}
	Let $W$ be a nonnegative function and $(\tilde{\gamma}_n)_{n \ge 0}$ be a $(\F_n)_{n \ge 0}$--predictable sequence of positive and finite random variables satisfying: there exist $\alpha >0$, $\beta > 0$, $n_0 \in \Nset$, such that  
	\begin{equation} \label{lem-condition-rappel-algo}
	    \forall n \ge n_0, \quad \frac{\E \bigl[ W(X_{n+1}) \,|\, \F_n \bigr] - W(X_n)}{ \tilde{\gamma}_{n+1} } \le - \alpha W(X_n) + \beta.
	\end{equation}
	Suppose $(\theta_n)_{n \ge 1}$ is a positive nonincreasing sequence such that $\E \bigl[ \sum_{n \ge 0} \theta_n \tilde{\gamma}_{n} \bigr]$ is finite, then
	\begin{equation} \label{lem-clef-controle-esp} 
		\E \Bigl[ \sum_{n \ge n_0 + 1} \theta_n \tilde{\gamma}_n W(X_{n-1}) \,|\, \F_{n_0} \Bigr] < + \infty.
	\end{equation}
	If, in addition, $\lim_n \theta_n = 0$ then 
	\begin{equation} \label{lem-clef-controle-vitesse}
		\lim_n \theta_n W(X_n) = 0 \quad a.s.
	\end{equation}
\end{prop}

The above proposition is related to Robbins-Siegmund's theorem (see Theorem 1.3.12 in \cite{duflo-rim} and the references therein).  

\begin{pf*}{Proof.}
	Let $\displaystyle R_n = \sum_{k = 0}^n \theta_k \tilde{\gamma}_k$ and $\displaystyle R_{\infty} = \lim_n R_n \in \L^1(\P)$.
	
	We consider the sequence $(Z_n)_{n \ge n_0}$ defined by  
	\begin{equation*}
		\forall n \ge n_0, \; Z_{n+1} = Z_n + \theta_{n+1} \bigl( W(X_{n+1}) - W(X_n) \bigr), \qquad Z_{n_0} = \theta_{n_0} W(X_{n_0}). 
	\end{equation*}	
	We first prove that for every $n \ge n_0$, $Z_n \ge 0$. Indeed, an Abel transform yields for every $n \ge n_0$, 
	\begin{align*}
	    Z_n &= \sum_{k = n_0}^{n-1} \theta_{k+1} \bigl( W(X_{k+1}) - W(X_k) \bigr) + \theta_{n_0} W(X_{n_0}), \\ 
	    &= \sum_{k = n_0}^{n-1} \bigl( \theta_k - \theta_{k + 1} \bigr) W(X_k) + \theta_n W(X_n).
	\end{align*}
	The sequence $(\theta_n)_{n \ge 0}$ is nonincreasing and the function $W$ is nonnegative, then $(Z_n)_{n \ge n_0}$ is positive. 
	
	Let $(S_n)_{n \ge n_0}$ denote the process defined for every $n \ge n_0$ by 
	\begin{equation*}
		S_n = Z_n + \alpha \sum_{k = n_0 + 1}^n \theta_k \tilde{\gamma}_k W(X_{k-1}) + \beta \bigl( \E \bigl[ R_{\infty} | \F_n \bigr] - R_n \bigr). 
	\end{equation*}
	Since $\E \bigl[ R_{\infty} \,|\, \F_n \bigr] - R_n \ge 0$, the sequence $(S_n)_{n \ge n_0}$ is nonnegative. Moreover, as $W$ satisfies \eqref{lem-condition-rappel-algo} we have 
	\begin{equation*}
		\forall n \ge n_0, \quad \E \bigl[ Z_{n+1} \,|\, \F_n \bigr] \le Z_n - \alpha \theta_{n+1} \tilde{\gamma}_{n+1} W(X_n) + \beta \theta_{n+1} \tilde{\gamma}_{n+1}.
	\end{equation*}
	Then it follows from this and from the $\F_n$--measurability of $\tilde{\gamma}_{n+1}$ that 
	\begin{equation*}
		\forall n \ge n_0, \quad \E \bigl[ S_{n+1} \,|\, \F_n \bigr] \le S_n.
	\end{equation*}
	Thus $(S_n)_{n \ge 0}$ converges $a.s.$ to a nonnegative finite random variable $S_{\infty}$, and we have
	\begin{equation*}
		\E \Bigl[ \sum_{n \ge n_0 + 1} \theta_n \tilde{\gamma}_n W(X_{n-1}) \,|\, \F_{n_0} \Bigr] < + \infty.
	\end{equation*}

	From the almost sure convergence of $(S_n)_{n \ge n_0}$, we also deduce the almost sure convergence of the series 
	\begin{equation*}
		\displaystyle \sum_{n \ge 1} \theta_n \bigl( W(X_n) - W(X_{n-1}) \bigr).
	\end{equation*}
	Since $(\theta_n)_{n \ge 0}$ is nonincreasing and converges to 0, Kronecker's lemma implies the almost sure convergence of $\bigl( \theta_n W(X_n) \bigr)_{n \ge 0}$ to $0$. 
\qed \end{pf*}

\begin{rem} A substitute for the $L^p$-boundedness of $(V(X_n))_{n \ge 0}$ has been already found in Lemma~4 of \cite{lamberton-pages-03} but does not apply in our case. Indeed, we have no information on the expectation of the random variable $\tilde{\gamma}_n$.
\end{rem}


The following Proposition is a fundamental consequence of \eqref{lem-clef-controle-vitesse} and says that for a ``good choice'' of the process $(\chi_n)_{n \ge 1}$, we can choose a sequence $(\gamma_n)_{n \ge 0}$ such that the random step $\tilde{\gamma}$ becomes deterministic after an almost surely finite time. The existence of this finite time $n_1$ is a significant property of our scheme.
\begin{prop} \label{prop-fondamentale}
	Let $f:[1, +\infty) \rightarrow (0, +\infty)$ be a decreasing one-to-one continuous function with $\displaystyle \lim_{x \rightarrow +\infty} f(x) = 0$ and $\displaystyle \lim_{x \rightarrow 1} f(x) = +\infty$, and $W$ a positive function with values in $[1, +\infty)$ satisfying \eqref{lem-condition-rappel-algo}. 
	If 
	\begin{equation*}
	    \chi_n \ge (f \circ W)(X_n), 
	\end{equation*}
	and if $(\gamma_n)_{n \ge 0}$ is subject to the condition 
	\begin{equation*}
	    \sum_{n \ge 1} \frac{\gamma_n}{f^{-1}(\gamma_n)} < + \infty,
	\end{equation*}
	where $f^{-1}$ is the inverse of $f$, then there is an almost surely finite random variable $n_1$ such that $\tilde{\gamma}_n = \gamma_n$ for every $n \ge n_1$.
\end{prop}
	
\begin{pf*}{Proof.}
	Let $(\theta_n)_{n \ge 0}$ the sequence defined by 
	\begin{equation}
		\forall n \ge 0, \quad \theta_n = \frac{1}{f^{-1}(\gamma_n)}.
	\end{equation}
	Since $(\gamma_n)_{n \ge 0}$ is nonincreasing and $1/f^{-1}$ is increasing on $\Rset_+$, then $(\theta_n)_{n \ge 0}$ is nondecreasing and we have $\sum_n \theta_n \tilde{\gamma}_n \le \sum_n \frac{\gamma_n}{f^{-1}(\gamma_n)} < +\infty$. Moreover $\lim_n \theta_n = 1 / f^{-1}(0) = 0$ so that, by Proposition \ref{lem-clef}, we have 
	\begin{equation*}
		\lim_n \theta_{n-1} W(X_{n-1}) = 0 \quad a.s. 
	\end{equation*}
	Hence there is an $a.s.$ finite random variable $n_1$ such that 
	\begin{equation*}
		\forall n \ge n_1, \quad W(X_{n-1}) < f^{-1}(\gamma_{n-1}) \le f^{-1}(\gamma_n) \quad a.s.
	\end{equation*}
	By the lower bound on $\chi_{n}$ and the monotony of $f$, we have 
	\begin{equation*}
		\forall n \ge n_1, \quad \chi_{n-1} > \gamma_n \quad a.s.
	\end{equation*}
	which completes the proof. 
\qed \end{pf*}
\section{Convergence of empirical measures}

In this section, we give conditions for the almost sure tightness of $\bigl( \nu^{\eta}_n \bigr)_{n \ge 1}$ and the weak convergence to an invariant distribution of \eqref{intro-diffusion}. To this end, we assume the existence of an almost surely finite random variable $n_1$ such that for every $n \ge n_1$, $\tilde{\gamma}_n = \gamma_n$. In practice, this means that $\tilde{\gamma}_n$ is defined by \eqref{definition-gamma-tilde} with $\gamma_n$ and $\chi_n$ satisfying conditions of Proposition~\ref{prop-fondamentale}. Under this assumption the proofs are very close to those in \cite{lamberton-pages-02} and \cite{lamberton-pages-03}. 

From now on we make the assumption: 
\begin{assumption}
	The deterministic sequences $(\gamma_n)_{n \ge 0}$ and $(\eta_n)_{n \ge 0}$ satisfy
	\begin{equation*}
		\lim_n \frac{1}{H_n} \sum_{k = 1}^n \Bigl| \Delta \frac{\eta_k}{\gamma_k} \Bigr| = 0,  
	\end{equation*}
	\begin{equation} \label{hypothese-theta1}
		\biggl( \frac{1}{\gamma_n H_n} \Bigl(\Delta \frac{\eta_n}{\gamma_n} \Bigr)_+ \biggr) \text{ is nonincreasing and } \sum_n \frac{1}{H_n} \Bigl(\Delta \frac{\eta_n}{\gamma_n} \Bigr)_+ < + \infty, 
	\end{equation}
	and there exists $s \in (1, 2]$ such that 
	\begin{equation} \label{hypothese-theta2}
		\biggl( \frac{1}{\gamma_n} \Bigl( \frac{\eta_n}{H_n \sqrt{\gamma_n}} \Bigr)^s \biggr) \text{ is nonincreasing and } \sum_n \Bigl( \frac{\eta_n}{H_n \sqrt{\gamma_n}} \Bigr)^s < +\infty.
	\end{equation}
\end{assumption}

\begin{rem} In practice the above conditions on $(\gamma_n)_{n \ge 0}$ and $(\eta_n)_{n \ge 1}$ are not restrictive. Setting $\gamma_n = n^{-p}$, $0 < p \le 1$ and $\eta_n = n^{-q}$, $q \le 1$, Assumption~2 is satisfied if and only if 
\begin{equation*}
	(p, q) \in \Bigl(0, \frac{2 (s-1)}{s} \Bigr) \times (-\infty, 1] \cup \Bigl\{ \Bigl(\frac{2(s-1)}{s}, 1 \Bigr) \Bigr\}.
\end{equation*}
\end{rem}

\subsection{A.s. tightness of empirical measures} 
We begin with proving the almost surely tightness of $\bigl(\nu^{\eta}_n \bigr)_{n \ge 1}$ when we have a control on the scheme $(X_n)_{n \ge 0}$. 
\begin{thm} \label{thm-tension}
	We assume that $(\gamma_n)_{n \ge 0}$ and $(\eta_n)_{n \ge 0}$ satisfy Assumption 2 and that there exists an almost surely finite random variable $n_1$ such that $\tilde{\gamma}_n = \gamma_n$ for every $n \ge n_1$.
	Suppose that $W$ is a positive function satisfying \eqref{lem-condition-rappel-algo} and
	\begin{equation*}
		\sum_{n \ge 1} \Bigl( \frac{\eta_n}{H_n \gamma_n} \Bigr)^s \E \Bigl[ \bigl| W(X_n) - \E \bigl[ W(X_n) \,|\, \F_{n-1} \bigr] \bigl|^s \,|\, \F_{n-1} \Bigr] < + \infty \quad a.s.,
	\end{equation*}
	for some $s \in (1,2]$. Then, 
	\begin{equation*} 
		\sup_{n \ge 1} \nu^{\eta}_n (W) < +\infty \quad a.s.
	\end{equation*}
\end{thm}

\begin{pf*}{Proof.}
	Since $W$ satisfies \eqref{lem-condition-rappel-algo}, there exist $\alpha > 0$, $\beta > 0$ and $n_0 \ge 0$ such that 
	\begin{equation*}
		\forall n \ge n_0, \quad W(X_n) \le \frac{W(X_n) - \E \bigl[ W(X_{n+1}) \,|\, \F_n \bigr]}{\tilde{\gamma}_{n+1} \alpha} + \frac{\beta}{\alpha}.
	\end{equation*}
	Hence, for every $n \ge n_0 \vee n_1 + 1$, 
	\begin{multline*}
		\frac{1}{H_n} \sum_{k = n_0 \vee n_1 + 1}^n \eta_k W(X_{k-1}) \\ \le \frac{1}{H_n} \sum_{k = n_0 \vee n_1+ 1}^n \frac{\eta_k}{\gamma_k \alpha} \bigl( W(X_{k-1}) - \E \bigl[ W(X_k) \,|\, \F_{k-1} \bigr] \bigr) + \frac{\beta}{\alpha}.
	\end{multline*}
	It suffices to prove that 
	\begin{equation*}
		\sup_{n \ge 1} \frac{1}{H_n} \sum_{k = 1}^n \frac{\eta_k}{\gamma_k} \bigl( W(X_{k-1}) - \E \bigl[ W(X_k) \,|\, \F_{k-1} \bigr] \bigr) < + \infty \quad a.s.
	\end{equation*}
	
	An Abel transform, setting $\eta_0 = 0$, yields  
	\begin{align*}
		\frac{1}{H_n} \sum_{k = 1}^n \frac{\eta_k}{\gamma_k} \bigl( W(X_{k-1}) - W(X_k) \bigr) &= \frac{1}{H_n} \sum_{k=1}^n \Bigl( \Delta \frac{\eta_k}{\gamma_k} \Bigr) W(X_{k-1}) - \frac{\eta_n}{\gamma_n} W(X_n) \\
		&\le \frac{1}{H_n} \sum_{k = 1}^n \Bigl( \Delta \frac{\eta_k}{\gamma_k} \Bigr)_+ W(X_{k-1}),
	\end{align*}
	and it follows from condition \eqref{hypothese-theta1} and Proposition \ref{lem-clef} applied with $\theta_n = \frac{1}{\gamma_n H_n} \bigl(\Delta \frac{\eta_n}{\gamma_n}\bigr)_+$ that 
	\begin{equation*}
		\sum_{n \ge 1} \frac{1}{H_n} \Bigl( \Delta \frac{\eta_n}{\gamma_n} \Bigr)_+ W(X_{n-1}) < + \infty \quad a.s.
	\end{equation*}
	Applying Kronecker's lemma, we get 
	\begin{equation*}
		\limsup_n \frac{1}{H_n} \sum_{k = 1}^n \frac{\eta_k}{\gamma_k} \bigl( W(X_{k-1}) - W(X_k) \bigr) \le 0 \quad a.s.
	\end{equation*}
	 
	It remains to prove that
	\begin{equation*}
		\sup_{n \ge 1} \frac{1}{H_n} \sum_{k = 1}^n \frac{\eta_k}{\gamma_k} \bigl( W(X_k) - \E \bigl[ W(X_k) \,|\, \F_{k-1} \bigr] \bigr) < + \infty \quad a.s.
	\end{equation*}
	We introduce the martingale $(M_n)_{n \ge 1}$ defined by 
	\begin{equation*}
		\forall n \ge 1, \quad M_n := \sum_{k = 1}^n \frac{\eta_k}{H_k \gamma_k} \Bigl( W(X_k) - \E \bigl[ W(X_k) \, |\, \F_{k - 1} \bigr] \Bigr), \quad M_0  := 0. 
	\end{equation*}
	By the Chow theorem (see \cite{hall-heyde}), the $a.s.$ convergence of $(M_n)_{n \ge 1}$ will follow from the $a.s.$ convergence of  
	\begin{equation*}
		\sum_{n \ge 1} \Bigl( \frac{\eta_n}{H_n \gamma_n} \Bigr)^s \E \Bigl[ \bigl| W(X_n) - \E \bigl[ W(X_n) \,|\, \F_{n-1} \bigr] \bigl|^s \,|\, \F_{n-1} \Bigr] < + \infty \quad a.s.
	\end{equation*}
	and the Kronecker lemma completes the proof.
\qed \end{pf*}

\subsection{Identification of the limit} \label{subsection:identif-limite}
We now prove that any weak limit of $(\nu^{\eta}_n)_{n \ge 0}$ is an invariant distribution for the diffusion \eqref{intro-diffusion}. We use the same method as Lamberton and Pag\`{e}s. By the Echeverria-Weiss theorem (see \cite{hall-heyde}) it suffices to prove that $\lim_n \nu^{\eta}_n (\A f) = 0$ for any twice continuously differentiable function $f$ with compact support. The existence of the Lyapounov function $V$ implies the regularity of the process $(x_t)_{t \ge 0}$ so that the Echeverria-Weiss applies --although $b$ has not sublinear growth--.

\begin{prop} \label{thm-identif}
	Suppose that there exists $K > 0$ such that 
	\begin{equation} \label{hypo-acc-algo}
		\forall n \ge 1, \quad \abs{X_n - X_{n-1}} \le K \sqrt{\gamma_n} \sqrt{V}(X_{n-1}) \bigl( 1 + \abs{U_n} \bigr),
	\end{equation} 
	and a function $W$ satisfying \eqref{lem-condition-rappel-algo}, $V = o(W)$ and $\sup_{n \ge 1} \nu^\eta_n(W) < +\infty$.
	Then for every twice continuously differentiable function $f$ with compact support, 
	\begin{equation*}
		\lim \nu_n^{\eta}(\A f) = 0 \quad a.s.
	\end{equation*}
\end{prop}


The following Lemma is useful to prove this Proposition.
\begin{lem} Under the assumptions of Proposition~\ref{thm-identif}, for every bounded Lipschitz continuous function $f:\Rset^d \rightarrow \Rset$ we have 
\begin{equation*}
	\lim_n \frac{1}{H_n} \sum_{k = 1}^n \eta_k \frac{\E \bigl[ f(X_k) \,|\, \F_{k-1} \bigr] - f(X_{k-1})}{\gamma_k} = 0 \quad a.s.
\end{equation*}
\end{lem}

\begin{pf*}{Proof.}
First, we prove that 
\begin{equation} \label{idlim-lem-part2}
	\lim_n \frac{1}{H_n} \sum_{k = 1}^n \eta_k \frac{\E \bigl[ f(X_k) \,|\, \F_{k-1} \bigr] - f(X_k)}{\gamma_k} = 0 \quad a.s.
\end{equation}
To this end, we introduce the martingale $(M_n)_{n \ge 0}$ defined by  
\begin{equation*}
	\forall n \ge 1, \quad M_n = \sum_{k = 1}^n \frac{\eta_k}{H_k \gamma_k} \bigl( \E \bigl[ f(X_k) \,|\, \F_{k-1} \bigr] - f(X_k) \bigr) \quad \text{and} \quad M_0 = 0.
\end{equation*}
We have 
\begin{multline*}
	\bigl( \E \bigl[ f(X_k) \,|\, \F_{k-1} \bigr] - f(X_k) \bigr)^2  = \\ \bigl( \E \bigl[ f(X_k) - f(X_{k-1}) \,|\, \F_{k-1} \bigr] + f(X_{k-1}) - f(X_k) \bigr)^2 ,   
\end{multline*}
and using $(a+b)^2 \le 2(a^2 + b^2)$, Jensen's inequality and $f$ Lipschitz we get 
\begin{align}
	\psca{M}_{\infty} &\le 4 \sum_{n \ge 1} \Bigl( \frac{\eta_n}{H_n \gamma_n} \Bigr)^2 \E \Bigl[ \bigl( f(X_n) - f(X_{n-1}) \bigr)^2 \,|\, \F_{n-1} \Bigr], \notag \\
	& \le C \sum_{n \ge 1} \Bigl( \frac{\eta_n}{H_n \sqrt{\gamma_n}} \Bigr)^2 V(X_{n-1}). \label{pf-thm-identif-1}
\end{align}
Since $s \le 2$ and $(\gamma_n)_{n}$ is nonincreasing, the sequence $(\theta_n)_{n} = \bigl(\frac{1}{\gamma_n} \bigl( \frac{\eta_n}{H_n \sqrt{\gamma_n}} \bigr)^2\bigr)_{n \ge 1}$ is nonincreasing and by \eqref{hypothese-theta2} we have $\sum_{n \ge 1} \theta_n \gamma_n < +\infty$. Moreover $W$ satisfies $\eqref{lem-condition-rappel-algo}$ and $\tilde{\gamma}_n = \gamma_n$ for every $n \ge n_1$, then the Proposition~\ref{lem-clef} applied with $\gamma_n$ and $\theta_n$ yields 
\begin{equation} \label{pf-thm-identif-2}
	\sum_{n \ge n_0 \vee n_1 + 1} \theta_n \gamma_n W(X_{n-1}) < +\infty \quad a.s.
\end{equation}
By \eqref{pf-thm-identif-1}, $V = o(W)$ and \eqref{pf-thm-identif-2} we obtain the almost sure convergence of the increasing process of $(M_n)_{n \ge 1}$. Thus the martingale converges almost surely. The Kronecker lemma gives \eqref{idlim-lem-part2}.

Finally we prove that  
\begin{equation*}
		\lim_n \frac{1}{H_n} \sum_{k = 1}^n \frac{\eta_k}{\gamma_k} \bigl( f(X_k) - f(X_{k-1}) \bigr) = 0 \quad a.s., 
\end{equation*}
using Abel's transform, the boundedness of $f$ and $\lim_n \frac{1}{H_n} \sum_{k=1}^n \abs{\Delta \frac{\eta_k}{\gamma_k}} = 0$.
\qed \end{pf*}

\begin{pf*}{Proof of Proposition~\ref{thm-identif}.}
	By Taylor's formula applied to $f$ between $X_{n-1}$ and $X_n$ we have 
	\begin{multline} \label{taylor}
		\E \bigl[ f(X_n) \,|\, \F_{n-1} \bigr] - f(X_{n-1}) = \tilde{\gamma}_n \A f(X_{n-1}) + \frac{1}{2} \tilde{\gamma}^2_n D^2 f(X_{k-1}) \cdot b(X_{n-1})^{\otimes 2} \\ + \E \Bigl[ R_2(X_{n-1}, X_n) \,|\, \F_{n-1} \Bigr], 
	\end{multline}
	with $R_2(x, y) = f(y) - f(x) - \psca{\nabla f(x), y - x} - \frac{1}{2} D^2f(x) \cdot (y-x)^{\otimes 2}$. On the one hand, by the above lemma 
	\begin{equation} \label{limite1-nulle}
		\lim_n \frac{1}{H_n} \sum_{k = 1}^n \frac{\eta_k}{\tilde{\gamma}_k} \Bigl( \E \bigl[ f(X_k) \,|\, \F_{k-1} \bigr] - f(X_{k-1}) \Bigr) = 0 \quad a.s.
	\end{equation}	
	and on the other hand we have, using that $\tilde{\gamma}_n \le \gamma_n$ and that $D^2 f$ has compact support,  
	\begin{equation*}
		\Bigl| \frac{1}{H_n} \sum_{k = 1}^n \eta_k \tilde{\gamma}_k D^2f(X_{k-1}) \cdot b^{\otimes 2}(X_{k-1}) \Bigr| \le \normsup{D^2 f \cdot b^{\otimes2}} \frac{1}{H_n} \sum_{k = 1}^n \eta_k \gamma_k 
	\end{equation*}
	and since $(\gamma_n)_{n \ge 0}$ is decreasing to 0, we obtain  
	\begin{equation} \label{limite2-nulle}
		\lim_n \frac{1}{H_n} \sum_{k = 1}^n \eta_k \tilde{\gamma}_k D^2f(X_{k-1}) \cdot b^{\otimes 2}(X_{k-1}) = 0 \quad a.s.  
	\end{equation}
	From \eqref{taylor}, \eqref{limite1-nulle} and \eqref{limite2-nulle}, it follows that 
	\begin{equation} \label{identif-pre-result} 
		\lim_n \Bigl( \nu^{\eta}_n(\A f) + \frac{1}{H_n} \sum_{k = 1}^n \frac{\eta_k}{\gamma_k} \E \Bigl[ R_2(X_{k-1}, X_k) \,|\, \F_{k-1} \Bigr] \Bigr) = 0 \quad a.s. 
	\end{equation}
	
	We introduce the continuous bounded function  
	\begin{equation*}
		r_2(x, \delta) = \frac{1}{2} \sup_{z \in \Rset^d, \abs{z-x} < \delta} \abs{D^2f(z) - D^2f(x)},
	\end{equation*}
	which is nondecreasing in $\delta$ and satisfies 
	\begin{equation*}
		\forall (x, y)\in \Rset^{2d}, \quad \abs{R_2(x, y)} \le r_2(x, \abs{y-x}) \abs{y-x}^2.
	\end{equation*}
	Thus, from \eqref{identif-pre-result} and \eqref{hypo-acc-algo} it suffices to prove that 
	\begin{equation} \label{identif-a-prouver}
		\lim_n \frac{1}{H_n} \sum_{k = 1}^n \eta_k V(X_{k-1}) \E \Bigl[ r_2(X_{k-1}, \delta_k(U_k)) \bigl(1 + \abs{U_k}^2 \bigr) \,|\, \F_{k-1} \Bigr] = 0 \quad a.s. 
	\end{equation}
	where $\delta_k(u) = \abs{\tilde{\gamma}_k b(X_{k-1}) + \sqrt{\tilde{\gamma}_k} \sigma(X_{k-1}) u}$.

	Let $u \in \Rset^m$ and $A > 0$. It is clear that the sequence $\bigl( \delta_k(u) 1_{\{\abs{X_{k-1}} \le A \}} \bigr)_{k \ge 1}$ converges almost surely to $0$ and thus 
	\begin{equation*}
		\lim_k r_2(X_{k-1}, \delta_k(u)) 1_{ \{ \abs{X_{k-1}} \le A \}} = 0 \quad a.s. 
	\end{equation*}
	Since 
	\begin{multline*}
		\frac{1}{H_n} \sum_{k=1}^n \eta_k V(X_{k-1}) r_2(X_{k-1}, \delta_k(u))(1+\abs{u}^2) 1_{\{\abs{X_{k-1}} \le A\}} \\ \le K \normsup{r_2} (1+\abs{u}^2), 
	\end{multline*}
	we have by the dominated convergence theorem 
	\begin{equation} \label{identif-a-prouver-sur-compact}
		\lim_n \frac{1}{H_n} \sum_{k=1}^n \eta_k V(X_{k-1}) \int_{\Rset^m} r_2(X_{k-1}, \delta_k(u))(1+\abs{u}^2) \mu(du) 1_{\{\abs{X_{k-1}} \le A\}} = 0 \quad a.s.
	\end{equation}
	where $\mu$ is the law of $U_1$.
	On the other hand, we have 
	\begin{multline*}
		\varlimsup_n \frac{1}{H_n} \sum_{k=1}^n \eta_k V(X_{k-1}) \E \bigl[ r_2(X_{k-1}, \delta_k(U_k))(1+\abs{U_k}^2) \,|\, \F_{k-1} \bigr] 1_{\{\abs{X_{k-1}} > A\}} 
		\\ \le (1+m) \normsup{r_2} \sup_{\abs{x} > A} \abs{V(x)/W(x)} \sup_n \nu^\eta_n(W),	
	\end{multline*}
	and letting $A \rightarrow +\infty$ and combining with \eqref{identif-a-prouver-sur-compact} we obtain \eqref{identif-a-prouver}.
\qed \end{pf*}
\section{Monotone and dissipative problems}
We now apply our results to monotone problems. In this section we assume that $V$ is essentially quadratic but the drift $b$ need not be globally Lipschitz. 
\begin{assumption} 
	The function $V$ satisfies 
	\begin{equation} \label{gradient-V-Lipschitz}
		\exists C_V > 0, \quad \abs{\nabla V}^2 \le C_V V, \quad \text{and} \quad \normsup{D^2 V} := \sup_{x \in \Rset^d} \norm{D^2V} < + \infty. 
	\end{equation}
\end{assumption}
Under this condition and Assumption~1, Mattingly, Stuart and Higham proved the geometric ergodicity of \eqref{intro-diffusion} when $\sigma$ is constant (see \cite{mattingly-02} for details). We assume that $\sigma$ satisfies  
\begin{equation} \label{sigma-negligeable}
	\exists C_{\sigma} > 0, \; \exists a \in (0, 1], \quad \Tr(\sigma \sigma^*) \le C_{\sigma} V^{1-a},
\end{equation}
so that condition \eqref{condition-sigma-V} is checked. 

We consider the Markov chain $(X_n)_{n \ge 0}$ built by the recursive procedure \eqref{algorithme} with the random step sequence $\tilde{\gamma}$ defined by 
\begin{equation}
	\tilde{\gamma}_0 = \gamma_0, \quad \forall n \ge 1,\; \tilde{\gamma}_n = g_n(X_{n-1})	
\end{equation}
with $g_n: \Rset^d \rightarrow \Rset_+^*$ defined by $g_n(x) = \gamma_n \wedge \chi(x)$ where $\chi$ is a function on $\Rset^d$ with value in $\Rset_+^*$.

The main result of this section is the following Theorem.
\begin{thm} \label{thm-essquad}
	Let $(\gamma_n)_{n \ge 0}$ and $(\eta_n)_{n \ge 0}$ satisfy Assumption $2$. Suppose that there exist $l > 1$ and $C_b \ge 1$ such that $\abs{b(x)}^2 \le C_b V^l(x)$ and that $\chi$ satisfies 
	\begin{equation} \label{thm-essquad-majo-chi}
		\forall x \in \Rset^d, \quad \zeta V^{-p}(x) \le \chi(x) \le \Bigl( \frac{2 \delta}{\normsup{D^2 V}} \Bigr) \frac{V(x)}{\abs{b(x)}^2 \vee 1},
	\end{equation}
	with $\zeta \le \displaystyle \frac{2 \delta}{C_b \normsup{D^2 V}}$, $p \ge l - 1$ and $\delta \in (0, \alpha)$.
	If $(\gamma_n)_{n \ge 0}$ satisfies 
	\begin{equation} \label{condition-gamma-essquad}
		\sum_{n \ge 1} \gamma_n \exp \bigl( -\lambda \zeta^{\frac{a}{p}} \gamma_n^{-\frac{a}{p}} \bigr) < +\infty, 
	\end{equation}
	and $\lambda_0 = \displaystyle \frac{2 \tau}{a C_\sigma \normsup{D^2 V}} \wedge \frac{2(\alpha - \delta)}{\kappa a C_V C_{\sigma}}$ then 
	\begin{itemize}
	\item there exists a finite random variable $n_1$ such that $\forall n \ge n_1$, $\tilde{\gamma}_n = \gamma_n$, 
	\item
	for every $\lambda < \frac{\lambda_0}{s}$ (whith $s$ as in \eqref{hypothese-theta2}), 
	\begin{equation*}
		\sup_n \nu^\eta_n \bigl( \exp (\lambda V^a) \bigr) < + \infty \quad a.s. 
	\end{equation*}
	and any weak limit of $\bigl( \nu^{\eta}_n \bigr)_{n \ge 1}$ is an invariant distribution for \eqref{intro-diffusion}.
	\end{itemize}
\end{thm}

\begin{rem} The condition \eqref{condition-gamma-essquad} is not restrictive. For example, it is satisfied if $\gamma_n \le \lambda^{\frac{p}{a}} \zeta \ln^{- \frac{p}{a}}(n)$ and $p > a$, or $\gamma_n \le \lambda^{\frac{p}{a}} \zeta \ln^{- (1+\varepsilon)}(n)$ for some $\varepsilon > 0$ and $p \le a$.
\end{rem}
To prove this Theorem, it suffices essentially to check that the function $W = \exp(\lambda V^a)$ satisfies condition~\eqref{lem-condition-rappel-algo} and the Assumptions of Theorem~\ref{thm-tension}. This is the aim of Lemmas~\ref{lem-rappel-expo-essquad} and \ref{lem-tension-essquad} respectively, which will be proved later.
\begin{lem} \label{lem-rappel-expo-essquad}
Assume that 
\begin{equation*} 
	\forall x \in \Rset^d, \quad \chi(x) \le \Bigl( \frac{2 \delta}{\normsup{D^2 V}} \Bigr) \frac{V(x)}{\abs{b(x)}^2 \vee 1},
\end{equation*}
where $\delta \in (0, \alpha)$. Let $\lambda_0$ be as in Theorem~\ref{thm-essquad}. Then for every $\lambda < \lambda_0$ there exists $\tilde{\alpha} > 0$, $\tilde{\beta} > 0$, and $n_0 \ge 0$ such that for every $n \ge n_0$, 
	\begin{equation*}
		\frac{ \E \Bigl[ \exp \bigl( \lambda V^a(X_{n+1}) \bigr) \,|\, \F_n \Bigr] - \exp \bigl( \lambda V^a(X_n) \bigr) }{ \tilde{\gamma}_{n+1} } \le - \tilde{\alpha} \exp \bigl( \lambda V^a(X_n) \bigr) + \tilde{\beta}.
	\end{equation*}
\end{lem}

\begin{lem} \label{lem-tension-essquad}
	Assume that there exists an almost sure finite random variable $n_1$ such that $\tilde{\gamma}_n = \gamma_n$ for every $n \ge n_1$. Let $\lambda_0$ be as in Theorem~\ref{thm-essquad}. Then for every $\displaystyle \lambda < \frac{\lambda_0}{s}$ the series  
	\begin{equation*}
		\sum_{n \ge 1} \Bigl( \frac{\eta_n}{H_n \gamma_n} \Bigr)^s \E \Bigl[ \bigl| \exp(\lambda V^a(X_n)) - \E \bigl[ \exp(\lambda V^a(X_n)) \,|\, \F_{n-1} \bigr] \bigl|^s \,|\, \F_{n-1} \Bigr] 
	\end{equation*}
	is almost surely finite.
\end{lem}

\begin{pf*}{Proof of Theorem~\ref{thm-essquad}.}
By Lemma~\ref{lem-rappel-expo-essquad} the function $W = \exp(\lambda V^a)$ satisfies condition \eqref{lem-condition-rappel-algo}.  
We consider the function $f$ defined on $[1, +\infty)$ by 
\begin{equation*}
	\forall x \ge 1, \quad f(x) = \zeta \lambda^{\frac{p}{a}} \ln^{-\frac{p}{a}}(x).
\end{equation*}
This is a decreasing one-to-one continuous function with $\displaystyle \lim_{x \rightarrow 1} f(x) = +\infty$ and $\displaystyle \lim_{x \rightarrow +\infty} f(x) = 0$. We check that $f \circ W = \zeta V^{-p}$ and by \eqref{thm-essquad-majo-chi} we have $\chi(X_n) \ge (f \circ W)(X_n)$.
The inverse of $f$ is the function $f^{-1}$ defined by 
\begin{equation}
	\forall x \in ]0, +\infty[, \quad f^{-1}(x) = \exp \bigl( \lambda \zeta^{\frac{a}{p}} x^{-\frac{a}{p}} \bigr),
\end{equation}
so that the condition \eqref{condition-gamma-essquad} on $\gamma_n$ is $\sum_{n \ge 0} \frac{\gamma_n}{f^{-1}(\gamma_n)} < +\infty$. The Proposition~\ref{prop-fondamentale} applied with $f$ and $W$ gives the existence of an almost surely finite random variable $n_1$ such that $\tilde{\gamma}_n = \gamma_n$ for every $n \ge n_1$.

By Lemma~\ref{lem-rappel-expo-essquad} the conditions of Theorem~\ref{thm-tension} are fulfilled and we have for every $\lambda < \lambda_0 / s$ 
\begin{equation*}
	\sup_{n \ge 1} \nu^\eta_n \bigl( \exp (\lambda V^a) \bigr) < + \infty \quad a.s. 
\end{equation*}
It remains to prove that any weak limit of $(\nu^\eta_n)_{n > 1}$ is an invariant distribution for the diffusion. By Subsection~\ref{subsection:identif-limite} and Proposition~\ref{thm-identif}, it suffices to check \eqref{hypo-acc-algo}. On the one hand, the definition of the algorithm yields  
\begin{equation*}
	\abs{X_n - X_{n-1}} \le \tilde{\gamma}_n \abs{b(X_{n-1})} + \sqrt{\tilde{\gamma}_n} \abs{\sigma(X_{n-1} U_n}.
\end{equation*}
On the other hand, we have $\tilde{\gamma}_n = g_n(X_{n-1})$ and 
\begin{equation} \label{monotone-ineg-fond}
	g_n(x) \abs{b(x)} \le \sqrt{g_n(x)} \sqrt{\chi(x)} \abs{b(x)} \le \sqrt{g_n(x)} \sqrt{\frac{2 \delta}{\normsup{D^2 V}}} \sqrt{V}(x),
\end{equation}
thus $\tilde{\gamma}_n \abs{b(X_{n-1})} \le \sqrt{\tilde{\gamma}_n} \sqrt{2 \delta/\normsup{D^2 V}} \sqrt{V}(X_{n-1})$. Since $\Tr(\sigma \sigma^*) \le C_\sigma V^{1-a}$, $a > 0$ and $V \ge 1$ then there exists $C > 0$ such that 
\begin{equation*}
	\abs{X_n - X_{n-1}} \le C \sqrt{\tilde{\gamma}_n} \sqrt{V}(X_{n-1}) ( 1 + \abs{U_n}),
\end{equation*}
and this, combined with $\tilde{\gamma}_n \le \gamma_n$, implies \eqref{hypo-acc-algo}.
\qed \end{pf*}

\begin{rem} The control of $\tilde{\gamma}_n \abs{b(X_{n-1})}$ given by \eqref{monotone-ineg-fond} is an important property of our scheme and will be often used throughout the section.  
\end{rem}

For the proof of Lemmas~\ref{lem-rappel-expo-essquad} and \ref{lem-tension-essquad} we will need the following consquence of conditions \eqref{bruit-sous-gaussien} and \eqref{bruit-moment-exponentiel} on $U_1$. 
\begin{lem} \label{lem-technique}
	There exists $K > 0$ such that $\forall \theta \in [0, \tau]$, $\forall h \in (0,1)$, $\forall v \in \Rset^d$,  
	\begin{equation*}
		\E \Bigl[ \exp \bigl(\sqrt{h} \psca{v, U_1} + h \theta \abs{U_1}^2 \bigr) \Bigr] \le \exp \Bigl( \frac{h}{1-h} \frac{\kappa}{2} \abs{v}^2 + K h \Bigr)
	\end{equation*}
\end{lem}

\begin{pf*}{Proof of Lemma~\ref{lem-technique}.}
	By H\"{o}lder's inequality, we have 
	\begin{multline*}
		\E \Bigl[ \exp \bigl(\sqrt{h} \psca{v, U_1} + h \theta \abs{U_1}^2 \bigr) \Bigr] \le \biggl( \E \Bigl[ \exp \Bigl(\frac{\sqrt{h}}{1-h} \psca{v, U_1} \Bigr) \Bigr] \biggr)^{1-h} \\
		\times \biggl( \E \Bigl[ \exp \bigl( \theta \abs{U_1}^2 \bigr) \Bigr] \biggr)^{h}.
	\end{multline*}
	Since the random variable $U_1$ is a generalized Gaussian, we have 
	\begin{equation*}
		\E \Bigl[ \exp \Bigl(\frac{\sqrt{h}}{1-h} \psca{v, U_1} \Bigr) \Bigr] \le \exp \Bigl( \frac{h}{(1-h)^2} \frac{\kappa}{2} \abs{v}^2 \Bigr). 
	\end{equation*}
	Since $U_1$ satisfies \eqref{bruit-moment-exponentiel} and $\theta \le \tau$, $\E \Bigl[ \exp \bigl( \theta \abs{U_1}^2 \bigr) \Bigr] < + \infty$. By setting 
	\begin{equation*}
		K = \log \Bigl( \E \Bigl[ \exp \bigl( \tau \abs{U_1}^2 \bigr) \Bigr] \Bigl), 
\end{equation*}
	the formula is established.
\qed \end{pf*}
\begin{pf*}{Proof of Lemma~\ref{lem-rappel-expo-essquad}.}

	We recall that $a \in ]0, 1]$. By concavity of the function $(x \mapsto x^a)$ we have 
	\begin{equation} \label{monotone-convexite}
		V^a(X_{n+1}) - V^a(X_n) \le a V^{a-1}(X_n) \bigl(V(X_{n+1}) - V(X_n) \bigr).
	\end{equation}
	The Taylor formula applied to $V$ between $X_n$ and $X_{n+1}$ yields
	\begin{equation} \label{taylor-V}
		V(X_{n+1}) \le V(X_n) + \psca{\nabla V(X_n), \Delta X_{n+1}} + \frac{1}{2} \rho \abs{\Delta X_{n+1}}^2, 
	\end{equation}
	where $\rho = \normsup{D^2 V}$.
	
	From the stability condition \eqref{condition-rappel} and $\Delta X_{n+1} = \tilde{\gamma}_{n+1} b(X_n) + \sqrt{\tilde{\gamma}_{n+1}} \sigma(X_n) U_{n+1}$, there exist $\alpha > 0$ and $\beta > 0$ such that
	\begin{equation} \label{taylor-V-majo-ordre1}
		\psca{\nabla V(X_n), \Delta X_{n+1}} \le - \alpha \tilde{\gamma}_{n+1} V(X_n) + \beta \tilde{\gamma}_{n+1} + \Lambda^{(1)}_{n+1}(X_n, U_{n+1}),
	\end{equation}
with the notation
	\begin{equation*}
		\Lambda_{n+1}^{(1)}(x, u) = \sqrt{g_{n+1}(x)} \psca{\nabla V(x), \sigma(x) u}.
	\end{equation*}
	
	On the other hand we write 
	\begin{equation*}
		\Lambda_{n+1}^{(2)}(x, u) = 2 \bigl( g_{n+1}(x) \bigr)^{3/2} \psca{b(x), \sigma(x) u} \\+ g_{n+1}(x) \Tr(\sigma \sigma^*)(x) \abs{u}^2,
	\end{equation*}
	so that 
	\begin{equation*}
		\abs{\Delta X_{n+1}}^2 \le \tilde{\gamma}_{n+1}^2 \abs{b(X_n)}^2 + \Lambda_{n+1}^{(2)}(X_n, U_{n+1}). 
	\end{equation*}
	From \eqref{monotone-ineg-fond} it follows that
	\begin{equation} \label{taylor-V-majo-ordre2}
		\abs{\Delta X_{n+1}}^2 \le \frac{2 \delta}{\rho} \tilde{\gamma}_{n+1} V(X_n) + \Lambda_{n+1}^{(2)}(X_n, U_{n+1}). 
	\end{equation}
	
	Combining \eqref{monotone-convexite}, \eqref{taylor-V}, \eqref{taylor-V-majo-ordre1} and \eqref{taylor-V-majo-ordre2} gives 
	\begin{multline} \label{taylor-V-majo-finale}
		V^a(X_{n+1}) \le V^a(X_n) - a (\alpha - \delta) \tilde{\gamma}_{n+1} V^a(X_n) + a \beta \tilde{\gamma}_{n+1} \\ + a V^{a-1}(X_n) \bigl(\Lambda_{n+1}^{(1)}(X_n, U_{n+1}) + \frac{\rho}{2} \Lambda_{n+1}^{(2)}(X_n, U_{n+1}) \bigr).
	\end{multline}
	Let $\lambda \in (0, \lambda_0)$. From \eqref{taylor-V-majo-finale} we deduce that for every $n \ge 0$  
	\begin{multline} \label{pre-result}
		\E \Bigl[ \exp \bigl( \lambda V^a(X_{n+1}) \bigr) \,|\, \F_n \Bigr] \\ \le \exp \Bigl( \lambda V^a(X_n) - \lambda a (\alpha - \delta) \tilde{\gamma}_{n+1} V^a(X_n) + \lambda a \beta \tilde{\gamma}_{n+1} \Bigr) Y_{n+1}
	\end{multline}
	where 
	\begin{equation*}
		Y_{n+1} = \E \Bigl[ \exp \Bigl(\lambda a V^{a-1}(X_n) \bigl( \Lambda_{n+1}^{(1)}(X_n, U_{n+1}) + \frac{\rho}{2} \Lambda_{n+1}^{(2)}(X_n, U_{n+1}) \bigr) \Bigr) \,|\, \F_n \Bigr].
	\end{equation*}

	We fix $n$ such that $\gamma_n < 1$ which will ensure $g_n(x) < 1$ for every $x \in \Rset^d$, and define $\phi_n(x)$ by 
	\begin{equation*}
		\phi_n(x) = \E \Bigl[ \exp \Bigl(\lambda a V^{a-1}(x) \bigl( \Lambda_n^{(1)}(x, U_1) + \frac{\rho}{2} \Lambda_n^{(2)}(x, U_1) \bigr) \Bigr) \Bigr],
	\end{equation*}
	so that $Y_{n+1} = \phi_{n+1}(X_n)$.
	Setting $v = \lambda a V^{a-1}(x) \sigma^*(x) \bigl( \nabla V(x) + \rho g_n(x) b(x) \bigr)$ we obtain by the definition of $\Lambda^{(1)}_n$ and $\Lambda^{(2)}_n$,  
	\begin{equation*}
		\phi_n(x) = \E \Bigl[ \exp \Bigl( \sqrt{g_n(x)} \psca{v , U_1} + \frac{\lambda a \rho}{2} V^{a-1}(x) g_n(x) \Tr(\sigma \sigma^*)(x) \abs{U_1}^2 \Bigr) \Bigr].
	\end{equation*}
	Since $\Tr(\sigma \sigma^*) \le C_{\sigma} V^{1-a}$ and $\lambda < \lambda_0 \le \frac{2 \tau}{a \rho C_{\sigma}}$, we are able to apply Lemma~\ref{lem-technique} with $\theta = \frac{\lambda a \rho C_{\sigma}}{2}$ and $h = g_n(x)$ which yields the existence of $K > 0$ such that 
  \begin{equation} \label{phi-holder}
		\phi_n(x) \le \exp \Bigl( \frac{g_n(x)}{1-g_n(x)} \frac{\kappa}{2} \abs{v}^2 + K g_n(x) \Bigr).
	\end{equation}
	
	It remains to handle $\abs{v}^2$. From $\abs{\nabla V}^2 \le C_V V$ and $\psca{\nabla V, b} \le \beta$ we obtain  
	\begin{equation*}
		\abs{v}^2 \le \lambda^2 a^2 V^{2a-2}(x) \Tr(\sigma \sigma^*)(x) \bigl( C_V V(x) + 2 \rho \beta g_n(x) + \rho^2 g_n^2(x) \abs{b(x)}^2 \bigr)
	\end{equation*}
	From \eqref{sigma-negligeable}, \eqref{monotone-ineg-fond} and $V \ge 1$ we have  
	\begin{equation*}
		\abs{v}^2 \le \lambda^2 a^2 C_\sigma V^{a-1}(x) \bigl( C_V V(x) + C g_n(x) V(x) \bigr), 
	\end{equation*}
	with $C = 2 \rho (\delta + \beta)$. In the following inequalities, the letter $C$ is used to denote a positive constant. Since 
	\begin{equation*}
		\frac{g_n(x)}{1-g_n(x)} = g_n(x) + \frac{g_n^2(x)}{1-g_n(x)} \le g_n(x) + \frac{g_n^2(x)}{1-\gamma_{n_0}},
	\end{equation*}
	where $n_0 = \min \{ n; \gamma_n < 1 \}$, it follows that 
	\begin{equation} \label{majo-sigmav}
		\frac{g_n(x)}{1-g_n(x)} \abs{v}^2 \le \lambda^2 a^2 C_V C_{\sigma} g_n(x) V^a(x) + C g_n^2(x) V^a(x).
	\end{equation}
	Combining \eqref{majo-sigmav} with \eqref{phi-holder} yields 
	\begin{equation} \label{majo-phi}
		\phi_n(x) \le \exp \Bigl( \lambda^2 a^2 C_V C_{\sigma} \frac{\kappa}{2} g_n(x) V^a(x) + C g_n^2(x) V^a(x) + K g_n(x) \Bigr), 
	\end{equation}
	and this inequality is true for every $n \ge n_0$. 
	
	Let $\displaystyle \bar{\alpha} = a(\alpha - \delta) - \lambda \frac{\kappa a^2 C_V C_{\sigma}}{2}$. Since $\displaystyle \lambda < \lambda_0 \le \frac{2 (\alpha - \delta)}{\kappa a C_V C_{\sigma}}$, we have $\bar{\alpha} > 0$, and combining \eqref{majo-phi} with \eqref{pre-result} we get for every $n \ge n_0$, 
	\begin{multline*}
		\E \Bigl[ \exp \bigl( \lambda V^a(X_{n+1}) \bigr) \,|\, \F_n \Bigr] \le \exp \Bigl( \lambda V^a(X_n) - \lambda \bar{\alpha} \tilde{\gamma}_{n+1} V^a(X_n) \\ + C \tilde{\gamma}_{n+1}^2 V^a(X_n) + \bar{K} \tilde{\gamma}_{n+1} \Bigr), 
	\end{multline*}
	where $\bar{K} = \lambda a \beta + K$. Setting $n_0' = \min \{n; \gamma_n < \lambda \bar{\alpha} / (2C) \}$ and $\tilde{\alpha} = \bar{\alpha}/2$ we have for every $n \ge n_0 \vee n_0'$, 
	\begin{equation*}
		\E \Bigl[ \exp \bigl( \lambda V^a(X_{n+1}) \bigr) \,|\, \F_n \Bigr] \le \exp \Bigl( \lambda V^a(X_n) - \lambda \tilde{\alpha} \tilde{\gamma}_{n+1} V^a(X_n) + K \tilde{\gamma}_{n+1} \Bigr). 
	\end{equation*}
	With the notation $\tilde{\beta} = \tilde{\alpha} \exp (K / \tilde{\alpha})$ and $n_0'' = \min \{ n; \gamma_n < 1 / \tilde{\alpha} \}$ we have by convexity of the exponential function: for every $n > n_0 \vee n_0' \vee n_0''$, 
	\begin{equation}
		\E \Bigl[ \exp \bigl( \lambda V^a(X_{n+1}) \bigr) \,|\, \F_n \Bigr] \le (1 - \tilde{\alpha} \tilde{\gamma}_{n+1} ) \exp \bigl( \lambda V^a(X_n) \bigr) + \tilde{\beta} \tilde{\gamma}_{n+1},
	\end{equation}
	which is the desired conclusion.
\qed \end{pf*}
	
\begin{pf*}{Proof of Lemma~\ref{lem-tension-essquad}.}
	Let $W = \exp \bigl( \lambda V^a \bigr)$ with $\lambda < \lambda_0 / s$. As $s > 1$, we have by convexity 
	\begin{multline*}
		\bigl| W(X_n) - \E \bigl[ W(X_n) \,|\, \F_{n-1} \bigr] \bigl|^s \\ \le 2^{s-1} \Bigl( \bigl| W(X_n) - W(X_{n-1}) \bigr|^s + \bigl| \E \bigl[ W(X_n) - W(X_{n-1}) \,|\, \F_{n-1} \bigr] \bigr|^s \Bigr),  
	\end{multline*}
	and by Jensen's inequality  
	\begin{equation*}
		\bigl| \E \bigl[ W(X_n) - W(X_{n-1}) \,|\, \F_{n-1} \bigr] \bigr|^s \le \E \bigl[ \bigl| W(X_n) - W(X_{n-1}) \bigr|^s \,|\, \F_{n-1} \bigr]. 
	\end{equation*}
	Thus it suffices to prove 
	\begin{equation*}
		\sum_{n \ge 1} \Bigl( \frac{\eta_n}{H_n \gamma_n} \Bigr)^s \E \Bigl[ \bigl| W(X_n) - W(X_{n-1}) \bigr|^s \,|\, \F_{n-1} \bigr] < +\infty \quad a.s.
	\end{equation*}
	Taylor's formula applied to the convex function $\bigl( x \mapsto \exp(\lambda x) \bigl)$ between $V^a(X_n)$ and $V^a(X_{n-1})$ yields 
	\begin{equation} \label{accr-W}  
		\abs{W(X_n) - W(X_{n-1})} \le \lambda \bigl( \exp( \lambda V^a(X_n)) + \exp(\lambda V^a(X_{n-1})) \bigr) \abs{\Delta V^a(X_n)}, 
	\end{equation}
	where $\Delta V^a(X_n) = V^a(X_n) - V^a(X_{n-1})$. As in the proof of Theorem~\ref{thm-essquad}, from \eqref{monotone-ineg-fond} and \eqref{sigma-negligeable} we deduce that 
	\begin{equation*}
		\abs{X_n - X_{n-1}} \le C \sqrt{\gamma_n} \sqrt{V}(X_{n-1}) \bigl( 1 + \abs{U_n} \bigr), 
	\end{equation*}
	and by Lemma 2. (a) of \cite{lamberton-pages-02} we get  
    \begin{equation} \label{accr-Va}
		\abs{V^a(X_n) - V^a(X_{n-1})} \le C \sqrt{\gamma_n} \bigl( \sqrt{V}(X_{n-1}) \bigr)^{2a \vee 1} \Bigl( 1 + \abs{U_n}^{2a \vee 1} \Bigr).
	\end{equation}
	To simplify notation, set $A_{n-1} = \bigl( \sqrt{V}(X_{n-1}) \bigr)^{2a \vee 1}$ and $B_n = 1 + \abs{U_n}^{2a \vee 1}$. Plugging \eqref{accr-Va} into \eqref{accr-W} we obtain 
	\begin{equation*}
		\abs{W(X_n) - W(X_{n-1})} \le \lambda C \sqrt{\gamma_n} \bigl( W(X_n) + W(X_{n-1}) \bigr) A_{n-1} B_n,
	\end{equation*}
	and 
	\begin{multline} \label{esp-accr-W}
		\E \bigl[ \abs{\Delta W(X_n)}^s \,|\, \F_{n-1} \bigr] \le \lambda^s C^s \gamma_n^{s/2} \Bigl( W^s(X_{n-1}) A^s_{n-1} \E \bigl[ B^s_n \,|\, \F_{n-1} \bigr] \\
		+ A^s_{n-1} \E \bigl[ W^s(X_n) B^s_n \,|\, \F_{n-1} \bigr] \Bigr).
	\end{multline}
	By Young's inequality for any $\varepsilon > 0$ 
	\begin{equation*}
		\E \bigl[ W^s(X_n) B^s_n \,|\, \F_{n-1} \bigr] \le \frac{1}{1+\varepsilon} \E \bigl[ W^{s(1+\varepsilon)}(X_n) \,|\, \F_{n-1} \bigr] + \frac{\varepsilon}{1+\varepsilon} K_\varepsilon 
	\end{equation*}
	with $K_\varepsilon = \E \bigl[ B_n^{s(1+\varepsilon)/\varepsilon} \,|\, \F_{n-1} \bigr]$. We choose $\varepsilon$ such that $\lambda s (1+\varepsilon) < \lambda_0$ and we get from Lemma~\ref{lem-rappel-expo-essquad} that for every $n \ge n_0$, 
	\begin{align*}
		\E \bigl[ W^s(X_n) B^s_n \,|\, \F_{n-1} \bigr] &\le \frac{1}{1+\varepsilon} \E \bigl[ \exp \bigl( \lambda s (1+\varepsilon) V^a(X_n) \bigr) \,|\, \F_{n-1} \bigr] + \frac{\varepsilon}{1+\varepsilon} K_\varepsilon, \\
		&\le \frac{1}{1+\varepsilon} \Bigl( \exp \bigl(\lambda s (1+\varepsilon) V^a(X_{n-1}) \bigr) + \tilde{\beta}\tilde{\gamma}_n + K_\varepsilon \varepsilon \Bigr),
	\end{align*}
	Combining this with \eqref{esp-accr-W} we obtain 
	\begin{multline*}
		\E \bigl[ \abs{ \Delta W(X_n) }^s \,|\, \F_n \bigr] \le \lambda^s C^s \gamma_n^{s/2} \Bigl( W^s(X_{n-1}) A_{n-1}^s \E [ B_1^s ] \\ + \frac{1}{1+\varepsilon} W^{s(1+\varepsilon)}(X_{n-1}) A^s_{n-1} + C_\varepsilon A^s_{n-1} \Bigr).
	\end{multline*}
	Keeping in mind that $A_{n-1} = \bigl( \sqrt{V}(X_{n-1}) \bigr)^{2a \vee 1}$ one checks that there exists $\bar{\lambda} \in (\lambda s (1+\varepsilon), \lambda_0)$ and $\tilde{K} > 0$ such that 
	\begin{equation*}
		\E \Bigl[ \abs{W(X_n) - W(X_{n-1})}^s \,|\, \F_{n-1} \Bigr] \le \tilde{K} \gamma_n^{s/2} \exp \bigl( \bar{\lambda} V^a(X_{n-1}) \bigr).
	\end{equation*}
	Consequently, it remains to prove that 
	\begin{equation} \label{fin-preuve-lem-identif-monotone}
		\sum_{n \ge 1} \Bigl( \frac{\eta_n}{H_n \sqrt{\gamma_n}} \Bigr)^s \exp \bigl( \bar{\lambda} V^a(X_{n-1}) \bigr) < + \infty \quad a.s. 
	\end{equation}
	We consider the nonincreasing sequence $(\theta_n)_{n \ge 1} = \bigl(\frac{1}{\gamma}_n \bigl( \frac{\eta_n}{H_n \sqrt{\gamma_n}} \bigr)^s \bigr)_{n \ge 1}$. From \eqref{hypothese-theta1} and $\tilde{\gamma}_n = \gamma_n$ for every $n \ge n_1$, we can apply the Proposition~\ref{lem-clef} with $\gamma_n$, $\theta_n$ and $W = \exp(\bar{\lambda} V^a)$ which gives \eqref{fin-preuve-lem-identif-monotone}.
\qed \end{pf*}
\section{Dissipative Hamiltonian systems}
In this section, we consider a stochastic differential system of the type 
\begin{equation} \label{systeme-hamiltonien}
	\begin{cases} 
		dq_t = \partial_p H(q_t, p_t) dt, \\ 
		dp_t = - \partial_q H(q_t, p_t) dt - F(q_t, p_t) \partial_p H(q_t, p_t) dt + c(q_t, p_t) dW_t, 
	\end{cases}
\end{equation}
with $(q(0), p(0)) = (q_0, p_0) \in \Rset^{2d}$. 
The process $(W_t)_{t \ge 0}$ is a $d$-dimensional Brownian motion and $c$ is continuous on $\Rset^{2d}$ with values in the set of $d \times d$ matrices. The \emph{Hamiltonian} $H$ is of class $\C^2$ on $\Rset^{2d}$ with real values  and the function $F$ is of class $\C^2$ on $\Rset^{2d}$ with values in the set of $d \times d$ matrices. 

We write this system in the abstract form \eqref{intro-diffusion} where 
\begin{gather*}
	x = \begin{pmatrix} 
		q \\ 
		p
	\end{pmatrix} \in \Rset^{2d}, \quad 
	b(x) = \begin{pmatrix} 
		\partial_p H(x) \\
		- \partial_q H(x) - F(x) \partial_p H(x) 
	\end{pmatrix} 
	= \begin{pmatrix} 
		b_1(x) \\ 
		b_2(x)
	\end{pmatrix}, 
	\\
	\text{and} \quad 
	\sigma(x) = \begin{pmatrix}
		0 \\
		c(x) 
	\end{pmatrix}.
\end{gather*}
We recall that we work always under Assumption~1. The existence of a Lyapounov function $V$ is a natural hypothesis for dissipative Hamiltonian systems, and in many cases we can determine $V$ using the Hamiltonian $H$. For example, for the Langevin equation (equation \eqref{systeme-hamiltonien} with $F(x) = \gamma > 0$, $c(x) = c$ invertible and $H(q,p) = \frac{\abs{p}^2}{2}+g(q)$ with $g\in \C^{\infty}(\Rset^d; \Rset)$ a polynomial function growing at infinity like $\abs{q}^{2l}$, $l \ge 1$), the Lyapounov function is defined for every $x=(q,p)\in \Rset^{2d}$ by $V(x) = H(x) + \frac{\gamma}{2}\psca{p,q} + \frac{\gamma^2}{4} \abs{q}^2 + 1$ (see \cite{mattingly-02}).

In this section we assume that:  
\begin{assumption} The function $V$ satisfies  
	\begin{equation} 
		\exists C_V > 0, \; \abs{\partial_p V}^2 \le C_V V \quad \text{and} \quad
		\sup_{(q,p) \in \Rset^{2d}} \norm{\partial_{pp}^2 V(q, p)} < + \infty.
	\end{equation}
\end{assumption}
Typical Lyapounov function of Hamiltonian systems may have arbitrary polynomial growth with respect to $q$ but are essentially quadratic with respect to $p$.
So in such a framework the following assumption is natural 
\begin{equation*}
	\Tr \bigl( \sigma^* (\nabla V)^{\otimes 2} \sigma \bigr) = \Tr \bigl( c^* (\partial_p V)^{\otimes 2} c \bigr) \le C_V V \Tr(c^* c),
\end{equation*}
and 
\begin{equation*}
	\sup_{x \in \Rset^{2d}} \Tr(\sigma^* D^2 V \sigma)(x) = \sup_{x \in \Rset^{2d}} \Tr(c^* \partial_{pp}^2 V c)(x) \le \rho \sup_{x \in \Rset^{2d}} \Tr( c^* c)(x),
\end{equation*}
where $\rho = \sup_{(q,p) \in \Rset^{2d}} \abs{\partial_{pp}^2 V(q, p)}$. Thus condition \eqref{condition-sigma-V} about $\sigma$ is satisfied as soon as  
\begin{equation}
	\exists C_{\sigma} > 0, \; a \in ]0, 1], \quad \Tr(c c^*) = \Tr(\sigma \sigma^*) \le C_{\sigma} V^{1-a}.
\end{equation}

These assumptions are very weak and are satisfied by a large class of examples derived from perturbed Hamiltonian systems. For a general model for the Hamiltonian and many examples (essentially multidimensional oscillators) we refer to \cite{soize} (page 10 for hypothesis on the Hamiltonian). 
    	
Our scheme $(X_n)_{n \ge 0}$ is built applying the recursive procedure $\eqref{algorithme}$ with the random step sequence $\tilde{\gamma}$ defined by
\begin{equation}
	\tilde{\gamma}_0 = \gamma_0, \quad \forall n \ge 1, \; \tilde{\gamma}_n = g_n(X_{n-1})
\end{equation}
with $g_n(x) = \gamma_n \wedge \chi_n(x)$ where $\chi_n: \Rset^{2d} \rightarrow \Rset_+^*$. Throughout the section, we consider the following function $\psi_n(x): \Rset^{2d} \rightarrow \Rset$ defined by 
\begin{equation}
	\forall x \in \Rset^d, \quad \psi_n(x) = \sup_{\bar{q} \in (q, q+\gamma_n b_1(x))} \norm{D^2 V(\bar{q}, p)} \vee 1.
\end{equation}

We introduce the notation $X_n = (Q_n, P_n)$ so that  
\begin{equation} \label{system-hamilton}
	\begin{cases}
		Q_{n+1} = Q_n + \tilde{\gamma}_{n+1} b_1(X_n), \\
		P_{n+1} = P_n + \tilde{\gamma}_{n+1} b_2(X_n) + \sqrt{\tilde{\gamma}_{n+1}} c(X_n) U_{n+1},
	\end{cases}
\end{equation}
and $(Q_0, P_0) = (q_0, p_0)$. Remark that $(Q_n)_{n \ge 0}$ is $(\F_n)_{n \ge 0}$--predictable.

The principal result is the following Theorem.
\begin{thm} \label{thm-hamilton}
Let $(\gamma_n)_{n \ge 1}$ and $(\eta_n)_{n \ge 1}$ satisfy Assumption 2. Suppose that there exist $l > 1$ and $C_b \ge 1$ such that $\psi_n(x) \abs{b(x)}^2 \le C_b V^l(x)$ and that $\chi$ satisfies  
\begin{equation} \label{thm-hamilton-majo-chi}
	\forall x \in \Rset^d, \quad \zeta V^{-p}(x) \le \chi_n(x) \le \frac{2 \delta V(x)}{\psi_n(x) \bigl( \abs{b(x)}^2 \vee 1 \bigr)},
\end{equation}
with $\zeta \le \frac{2 \delta}{C_b}$, $p \ge l-1$ and $\delta \in (0, \alpha/4)$. If $(\gamma_n)_{n \ge 0}$ satisfies 
\begin{equation}
	\sum_{n \ge 1} \gamma_n \exp \bigl(- \lambda \zeta^{\frac{a}{p}} \gamma_n^{-\frac{a}{p}} \bigr) < +\infty,
\end{equation}
and $\lambda_0 = \displaystyle \frac{\tau}{a C_\sigma \normsup{\partial_{pp}^2 V}} \wedge \frac{2(\alpha - 4 \delta)}{\kappa a C_V C_{\sigma}}$ then
\begin{itemize}
\item there exists an $a.s.$ finite random variable $n_1$ such that $\forall n \ge n_1$, $\tilde{\gamma}_n = \gamma_n$,
\item for every $\lambda < \frac{\lambda_0}{s}$ (whith $s$ as in \eqref{hypothese-theta2}), 
\begin{equation*}
	\sup_n \nu^\eta_n \bigl( \exp (\lambda V^a) \bigr) < + \infty \quad a.s. 
\end{equation*}
and any weak limit of $\bigl( \nu^{\eta}_n \bigr)_{n \ge 1}$ is an invariant distribution for \eqref{intro-diffusion}.
\end{itemize}
\end{thm}

The proof of this Theorem is essentially the same as the proof of Theorem~\ref{thm-essquad}. Lemma~\ref{lem-rappel-exponentiel-hamilton} gives condition on $\chi$ so that $W = \exp(\lambda V^a)$ satisfy \eqref{lem-condition-rappel-algo} and Lemma~\ref{lem-tension-hamilton} allows to apply Theorem~\ref{thm-tension}.

\begin{lem} \label{lem-rappel-exponentiel-hamilton}
Assume that 
\begin{equation*}
	\forall x \in \Rset^{2d}, \quad \chi_n(x) \le \frac{2 \delta V(x)}{\psi_n(x) \bigl( \abs{b(x)}^2 \vee 1 \bigr)},
\end{equation*}
where $\delta \in (0, \alpha/4)$. Let $\lambda_0$ be as in Theorem~\ref{thm-hamilton}. Then for every $\lambda < \lambda_0$ there exist $\tilde{\alpha} > 0$, $\tilde{\beta} > 0$, and $n_0 \ge 0$ such that for every $n \ge n_0$, 
	\begin{equation*}
		\frac{ \E \Bigl[ \exp \bigl( \lambda V^a(X_{n+1}) \bigr) \,|\, \F_n \Bigr] - \exp \bigl( \lambda V^a(X_n) \bigr)}{\tilde{\gamma}_{n+1}} \le - \tilde{\alpha} \exp \bigl( \lambda V^a(X_n) \bigr) + \tilde{\beta}.
	\end{equation*}
\end{lem}

\begin{lem} \label{lem-tension-hamilton}
Assume that there exists an almost sure finite random variable $n_1$ such that $\tilde{\gamma}_n = \gamma_n$ for every $n \ge n+1$. Let $\lambda_0$ be as in Theorem~\ref{thm-hamilton}. Then for every $\displaystyle \lambda < \frac{\lambda_0}{s}$ the series 
	\begin{equation*}
		\sum_{n \ge 1} \Bigl( \frac{\eta_n}{H_n \gamma_n} \Bigr)^s \E \Bigl[ \bigl| \exp(\lambda V^a(X_n)) - \E \bigl[ \exp(\lambda V^a(X_n)) \,|\, \F_{n-1} \bigr] \bigl|^s \,|\, \F_{n-1} \Bigr]
	\end{equation*}
is almost surely finite.
\end{lem}

A first property of our scheme is given by the following Lemma. 
\begin{lem} \label{lem-hamilton-essentiel} If $\delta \le \alpha$ and 
\begin{equation} \label{lem-hamilton-essentiel-hypothese}
	\forall x \in \Rset^{2d}, \quad \chi_n(x) \le \frac{2 \delta V(x)}{\psi_n(x) \bigl(\abs{b(x)}^2 \vee 1 \bigr)},
\end{equation}
then for every $x = (q,p) \in \Rset^{2d}$,  
\begin{equation*} 
	V(q+g_n(x) b_1(x), p) \le \bigl( 1 + \sqrt{2 \delta C_V} \sqrt{g_n(x)} \bigr) V(x) + g_n(x) \beta.
\end{equation*}
\end{lem}

\begin{pf*}{Proof.} 
By the Taylor's formula we obtain   
\begin{multline} \label{hamilton-tech-1}
	V(q+g_n(x) b_1(x), p) - V(x) \le g_n(x) \psca{\partial_q V(x), b_1(x)} \\ + \frac{1}{2} \norm{\partial_{qq}^2 V(\bar{q}_n, p)} \abs{g_n(x) b_1(x)}^2, 
\end{multline}
where $\bar{q}_n \in (q, q+g_n(x) b_1(x))$. Since $g_n(x) = \gamma_n \wedge \chi_n(x)$ we have 
\begin{equation}
	g_n^2(x) \abs{b_1(x)}^2 \le g_n(x) \chi_n(x) \abs{b_1(x)}^2 \le g_n(x) \frac{2 \delta V(x)}{\psi_n(x)}.
\end{equation}
Moreover, $\norm{\partial_{qq}^2 V(\bar{q}_n, p)} \le \norm{D^2 V(\bar{q}_n, p)} \le \psi_n(x)$ and by \eqref{hamilton-tech-1} 
\begin{equation} \label{hamilton-tech-pre-result}
	V(q+g_n(x) b_1(x), p) - V(x) \le g_n(x) \psca{\partial_q V(x), b_1(x)} + g_n(x) \delta V(x). 
\end{equation}

Writing $\psca{\partial_q V(x), b_1(x)} = \psca{\nabla V, b}(x) - \psca{\partial_p V, b_2}(x)$ and using the stability condition \eqref{condition-rappel} we get 
\begin{multline*}
	V(q+g_n(x) b_1(x), p) - V(x) \le g_n(x) \bigl(\beta - \alpha V(x) - \psca{\partial_p V(x), b_2(x)} \bigr) \\ + g_n(x) \delta V(x). 
\end{multline*}
Since $\delta \le \alpha$ and $\abs{\partial_p V}^2 \le C_V V$ we have 
\begin{equation*}
	V(q+g_n(x) b_1(x), p) - V(x) \le g_n(x) \sqrt{C_V} \sqrt{V}(x) \abs{b_2(x)} + g_n(x) \beta,
\end{equation*}
and if $\abs{b_2(x)} \le \frac{\sqrt{2 \delta}}{\sqrt{g_n(x)}} \sqrt{V}(x)$ then 
\begin{equation*}
	V(q+g_n(x) b_1(x), p) - V(x) \le \sqrt{g_n(x)} \sqrt{2 \delta C_V} V(x) + g_n(x) \beta,
\end{equation*}
else $\sqrt{V}(x) < \frac{\sqrt{g_n(x)}}{\sqrt{2 \delta}} \abs{b_2(x)}$ and 
\begin{align*}
	V(q+g_n(x) b_1(x), p) - V(x) &\le g_n(x) \sqrt{g_n(x)} \sqrt{\frac{C_V}{2 \delta}} \abs{b_2(x)}^2 + g_n(x) \beta, \\
	&\le \sqrt{g_n(x)} \sqrt{\frac{C_V}{2 \delta}} \chi_n(x) \abs{b_2(x)}^2 + g_n(x) \beta.
\end{align*}
From \eqref{lem-hamilton-essentiel-hypothese} we deduce that 
\begin{equation*}
	V(q+g_n(x) b_1(x), p) - V(x) \le \sqrt{g_n(x)} \sqrt{2 \delta C_V} V(x) + g_n(x) \beta,
\end{equation*}
and the Lemma is proved.
\qed \end{pf*}

\begin{pf*}{Proof of Lemma~\ref{lem-rappel-exponentiel-hamilton}.}
	First, remark that by the concavity of the function $(x \mapsto x^a)$ we have 
	\begin{equation} \label{hamilton-passage-Va-V}
		V^a(X_{n+1}) - V^a(X_n) \le a V^{a-1}(X_n) \Bigl(V(X_{n+1}) - V(X_n) \Bigr),
	\end{equation}
	and that 
	\begin{equation*} 
		V(X_{n+1}) - V(X_n) = V(X_{n+1}) - V(Q_{n+1}, P_n) + V(Q_{n+1}, P_n) - V(X_n).	
	\end{equation*}
	
	In the proof of Lemma~\ref{lem-hamilton-essentiel} we proved \eqref{hamilton-tech-pre-result} i.e. 
	\begin{equation} \label{hamilton-accr-V-1}
		V(Q_{n+1}, P_n) - V(X_n) \le \tilde{\gamma}_{n+1} \psca{\partial_q V(X_n), b_1(X_n)} + \tilde{\gamma}_{n+1} \delta V(X_n). 
	\end{equation}
	
	We now study $V(X_{n+1}) - V(Q_{n+1}, P_n)$. We apply Taylor's formula to $V$ which gives 
	\begin{multline} \label{hamilton-accr-V-2-init} 
		V(X_{n+1}) - V(Q_{n+1}, P_n) \le \psca{ \partial_p V(Q_{n+1}, P_n), \Delta P_{n+1}} + \\ \frac{1}{2} \sup_{p \in (P_n, P_{n+1})} \norm{\partial_{pp}^2 V(Q_{n+1}, p)} \abs{\Delta P_{n+1}}^2.
	\end{multline}
	On the one hand, we have 
	\begin{align*}
		\psca{\partial_p V(Q_{n+1}, P_n), \Delta P_{n+1}} & = \tilde{\gamma}_{n+1} \psca{\partial_p V(X_n), b_2(X_n)} \\
		& \qquad + \tilde{\gamma}_{n+1} \psca{\partial_p V(Q_{n+1}, P_n) - \partial_p V(X_n), b_2(X_n)} \\ 
		& \qquad {} + \sqrt{\tilde{\gamma}_{n+1}} \psca{\partial_p V(Q_{n+1}, P_n), c(X_n) U_{n+1}},
	\end{align*}
	and using 
	\begin{align*}
		\abs{\partial_p V(Q_{n+1}, P_n) - \partial_p V(X_n)} & \le \norm{\partial_{qp}^2 V(q, P_n)} \abs{\Delta Q_{n+1}}, \\
		& \le \tilde{\gamma}_{n+1} \norm{D^2 V(q, P_n)} \abs{b_1(X_n)},
	\end{align*}
	where $q \in (Q_n, Q_{n+1})$ we obtain 
	\begin{align*}
		\psca{\partial_p V(Q_{n+1}, P_n), \Delta P_{n+1}} & \le \tilde{\gamma}_{n+1} \psca{\partial_p V(X_n), b_2(X_n)} \\
		& \qquad + \tilde{\gamma}_{n+1}^2 \norm{D^2 V(q, P_n)} \abs{b_1(X_n)} \abs{b_2(X_n)} \\ 
		& \qquad {} + \sqrt{\tilde{\gamma}_{n+1}} \psca{\partial_p V(Q_{n+1}, P_n), c(X_n) U_{n+1}}.
	\end{align*}
	Since $\tilde{\gamma}_{n+1} = \gamma_{n+1} \wedge \chi_{n+1}(X_n)$ we have 
	\begin{equation}
		\tilde{\gamma}_{n+1} \le \chi_{n+1}(X_n) \le \displaystyle \frac{2 \delta V(X_n)}{\displaystyle \sup_{q \in (Q_n, Q_{n+1})} \norm{D^2 V(q, P_n)} \bigl(\abs{b(X_n)}^2 \vee 1 \bigr)}
	\end{equation}
	and from $\abs{b_1(X_n)}\abs{b_2(X_n)} \le \frac{1}{2} \abs{b(X_n)}^2$ we obtain  
	\begin{multline} \label{hamilton-accr-V-2-part1} 
		\psca{\partial_p V(Q_{n+1}, P_n), \Delta P_{n+1}} \le \tilde{\gamma}_{n+1} \psca{\partial_p V(X_n), b_2(X_n)} + \tilde{\gamma}_{n+1} \delta V(X_n)\\ + \sqrt{\tilde{\gamma}_{n+1}} \psca{\partial_p V(Q_{n+1}, P_n), c(X_n) U_{n+1}}.
	\end{multline}	
	On the other hand, setting $\rho = \sup_{x\in \Rset^{2d}} \norm{\partial_{pp}^2 V(x)}$ and using \eqref{system-hamilton} we have  
	\begin{multline*}
		\sup_{p \in (P_n, P_{n+1})} \norm{\partial_{pp}^2 V(Q_{n+1}, p)}\abs{\Delta P_{n+1}}^2 \\ \le 2 \tilde{\gamma}_{n+1}^2 \sup_{p \in (P_n, P_{n+1})} \norm{D^2 V(Q_{n+1}, p)} \abs{b_2(X_n)}^2 \\ + 2 \tilde{\gamma}_{n+1} \rho \Tr(c c^*)(X_n) \abs{U_{n+1}}^2.
	\end{multline*}
	By the definition of $\tilde{\gamma}_{n+1}$ it follows that 
	\begin{multline} \label{hamilton-accr-V-2-part2}
		\sup_{p \in (P_n, P_{n+1})} \norm{\partial_{pp}^2 V(Q_n, p)}\abs{\Delta P_{n+1}}^2 \le 4 \tilde{\gamma}_{n+1} \delta V(X_n) \\ + 2 \tilde{\gamma}_{n+1} \rho \Tr(c c^*)(X_n) \abs{U_{n+1}}^2.
	\end{multline}
	Finally, combining \eqref{hamilton-accr-V-2-init} with \eqref{hamilton-accr-V-2-part1} and \eqref{hamilton-accr-V-2-part2} we obtain 
	\begin{multline} \label{hamilton-accr-V-2}
		V(X_{n+1}) - V(Q_{n+1}, P_n) \\ \le \tilde{\gamma}_{n+1} \psca{ \partial_p V(X_n), b_2(X_n)} + 3 \delta \tilde{\gamma}_{n+1} V(X_n) + Y_{n+1},
	\end{multline}
	where 
	\begin{equation*}
		Y_{n+1} = \sqrt{\tilde{\gamma}_{n+1}} \psca{\partial_p V(Q_{n+1}, P_n), c(X_n) U_{n+1}} + \rho \tilde{\gamma}_{n+1} \Tr(c c^*)(X_n) \abs{U_{n+1}}^2.
	\end{equation*}
	
	From \eqref{hamilton-accr-V-1} and \eqref{hamilton-accr-V-2} we have 
	\begin{equation*}
		V(X_{n+1}) - V(X_n) \le \tilde{\gamma}_{n+1} \psca{\nabla V, b}(X_n) + 4 \delta \tilde{\gamma}_{n+1} V(X_n) + Y_{n+1},
	\end{equation*}
	and by the stability condition \eqref{condition-rappel} there exists $\alpha > 0$ such that 
	\begin{align*}
		V(X_{n+1}) &\le V(X_n) - \alpha \tilde{\gamma}_{n+1} V(X_n) + 4 \delta \tilde{\gamma}_{n+1} V(X_n) + \beta \tilde{\gamma}_{n+1} + Y_{n+1}, \\
		&= V(X_n) - (\alpha - 4 \delta) \tilde{\gamma}_{n+1} V(X_n) + \beta \tilde{\gamma}_{n+1} + Y_{n+1}.
	\end{align*}
	By \eqref{hamilton-passage-Va-V} and $V \ge 1$ we get (using that $V^{a-1} \le 1$) 
	\begin{equation*}
		V^a(X_{n+1}) \le V^a(X_n) - a (\alpha - 4 \delta) \tilde{\gamma}_{n+1} V^a(X_n) + a \beta \tilde{\gamma}_{n+1} + a V^{a-1}(X_n) Y_{n+1}.
	\end{equation*}

	Let $\lambda \in (0, \lambda_0)$. Denoting 
	\begin{equation*}
		Z_{n+1} = \E \Bigl[ \exp \bigl( a \lambda V^{a-1}(X_n) Y_{n+1} \bigr) \,|\, \F_n \Bigr],	
	\end{equation*}
	we have 	
	\begin{multline} \label{hamilton-pre-result}
		\E \Bigl[ \exp \bigl( \lambda V^a(X_{n+1}) \bigr) \,|\, \F_n \Bigr] \\ \le \exp \bigl( \lambda V^a(X_n) - a \lambda (\alpha - 4 \delta) \tilde{\gamma}_{n+1} V^a(X_n) + a \lambda \beta \tilde{\gamma}_{n+1} \bigr) Z_{n+1}.
	\end{multline}

	It remains to study $Z_n$. We define $\phi_n(x)$ by 
	\begin{equation}
		\phi_n(x) = \E \Bigl[ \exp \bigl( \sqrt{g_n(x)} \psca{v, U_1} + \lambda a \rho g_n(x) V^{a-1}(x) \Tr(c c^*)(x) \abs{U_1}^2 \bigr) \Bigr]
	\end{equation}
	with $v = \lambda a V^{a-1}(x) c^*(x) \partial_p V \bigl(q + g_n(x) b_1(x), p \bigr)$, so that $Z_{n+1} = \phi_{n+1}(X_n)$. Using $\Tr(c c^*) \le C_{\sigma} V^{1-a}$ we have 
	\begin{equation*}
		\phi_n(x) \le \E \Bigl[ \exp \bigl( \sqrt{g_n(x)} \psca{v, U_1} + \lambda a \rho g_n(x) C_\sigma \abs{U_1}^2 \bigr) \Bigr].
	\end{equation*}
	Let $n_0 = \min \{n \ge 0; \gamma_n < 1 \}$. Since $\lambda < \lambda_0 \le \frac{\tau}{a \rho C_{\sigma}}$, we are able to apply lemma \ref{lem-technique} with $\theta = \lambda a \rho C_{\sigma}$ and $h = g_n(x)$ which gives the existence of $K > 0$ such that for every $n \ge n_0$  
	\begin{equation}
		\phi_n(x) \le \exp \Bigl( \frac{g_n(x)}{1-g_n(x)} \frac{\kappa}{2} \abs{v}^2 + K g_n(x) \Bigr).
	\end{equation}
	Moreover, it follows from  $\abs{\partial_p V} \le C_V V$ that  
	\begin{align*}
		\abs{v}^2 & \le \lambda^2 a^2 C_{\sigma} V^{a-1}(x) \abs{\partial_p V(q+g_n(x) b_1(x), p)}^2, \\
	    	& \le \lambda^2 a^2 C_{\sigma} C_V V^{a-1}(x) V(q+g_n(x) b_1(x), p),
	\end{align*}
	The Lemma~\ref{lem-hamilton-essentiel} gives 
	\begin{equation*}
		\abs{v}^2 \le \lambda^2 a^2 C_{\sigma} C_V V^{a}(x) + C \sqrt{g_n}(x) V^a(x) + g_n(x) \beta.
	\end{equation*}
	In the same manner as in the proof of Lemma~\ref{lem-rappel-expo-essquad} there exists $C > 0$ and $K > 0$ such that for every $n \ge n_0$ 
	\begin{equation*}
		\phi_n(x) \le \exp \Bigl( \frac{\kappa}{2} \lambda^2 a^2 C_{\sigma} C_V g_n(x) V^a(x) + C g_n^{3/2}(x) V^a(x) + K g_n(x) \Bigr)
	\end{equation*}
	then from \eqref{hamilton-pre-result} we have 
	\begin{multline*} 
		\E \Bigl[ \exp \bigl( \lambda V^a(X_{n+1}) \bigr) \,|\, \F_n \Bigr] \le \exp \bigl( \lambda V^a(X_n) \\ - a \lambda (\alpha - 4 \delta - \frac{\kappa}{2} \lambda a C_\sigma C_V) \tilde{\gamma}_{n+1} V^a(X_n) \\ + C \tilde{\gamma}_{n+1}^{3/2} V^a(X_n) + K \tilde{\gamma}_{n+1} \bigr).
	\end{multline*} 
	Let $\bar{\alpha} = a(\alpha - 4 \delta - \frac{\kappa}{2} \lambda a C_\sigma C_V)$. Since $\displaystyle \lambda < \lambda_0 \le \frac{2(\alpha - 4\delta)}{\kappa a C_{\sigma} C_V}$ we have $\bar{\alpha} > 0$. Setting $n_0' = \min\{n; \sqrt{\gamma_n} < a \lambda \bar{\alpha}/(2C) \}$ and $\tilde{\alpha} = (\lambda a \bar{\alpha})/2$ we have for every $n \ge n_0 \vee n_0'$  
	\begin{equation*} 
		\E \Bigl[ \exp \bigl( \lambda V^a(X_{n+1}) \bigr) \,|\, \F_n \Bigr] \le \exp \bigl( \lambda V^a(X_n) - \tilde{\alpha} \tilde{\gamma}_{n+1} V^a(X_n) + C \tilde{\gamma}_{n+1} \bigr),
	\end{equation*}
	which proves the lemma (by the convexity of the exponential).
\qed \end{pf*}
\begin{pf*}{Proof of Lemma~\ref{lem-tension-hamilton}.}
	Let $W = \exp ( \lambda V^a )$ with $\lambda < \lambda_0 / s$. 
	We recall that $Q_n$ is $\F_{n-1}$ measurable. Since $s > 1$, the convexity of $(x \mapsto x^s)$ implies that 
	\begin{multline*}
		\E \Bigl[ \bigl| W(X_n) - \E \bigl[ W(X_n) \,|\, \F_{n-1} \bigr] \bigl|^s \,|\, \F_{n-1} \Bigr] \\ \le 2^s \E \Bigl[ \abs{W(X_n) - W(Q_n, P_{n-1})}^s \,|\, \F_{n-1} \Bigr]. 
	\end{multline*}	
	First we prove that 
	\begin{equation}\label{ac-Va-amontrer}
		\abs{V^a(X_n) - V^a(Q_n, P_{n-1})} \le C \sqrt{\gamma_n} A_{n-1} B_n 
	\end{equation}
	with 
	\begin{align*}
		A_{n-1} &= \bigl(\sqrt{V}(Q_n, P_{n-1})\bigr)^{(2a-1)_+} \bigl(\sqrt{V}(X_{n-1})\bigr)^{(2a \vee 1)}, \\
		B_n &= 1+\abs{U_n}^{2a \vee 1}.
	\end{align*}
	Since $\Delta P_n = \tilde{\gamma}_n b_2(X_{n-1}) + \sqrt{\gamma_n} c(X_{n-1}) U_n$ and $\Tr(c c^*) \le C_\sigma V^{1-a}$, we have  
	\begin{equation*} 
		\abs{\Delta P_n} \le \tilde{\gamma}_n \abs{b_2(X_{n-1})} + \sqrt{\tilde{\gamma}_n} C_\sigma \bigl(\sqrt{V}(X_{n-1})\bigr)^{1-a} \abs{U_n}.
	\end{equation*}
	Using $\tilde{\gamma}_n \le \chi_n(X_{n-1}) \le \frac{2 \delta V(X_{n-1})}{\abs{b(X_{n-1})}^2 \vee 1}$ we obtain 
	\begin{align} \label{accr-Pn} 
		\abs{\Delta P_n} & \le \sqrt{\tilde{\gamma}_n} \sqrt{2 \delta} \sqrt{V}(X_{n-1}) + \sqrt{\tilde{\gamma}_n} C_\sigma \bigl(\sqrt{V}(X_{n-1})\bigr)^{1-a} \abs{U_n}, \notag \\
		& \le C \sqrt{\tilde{\gamma}_n} \sqrt{V}(X_{n-1}) \bigl(1+\abs{U_n} \bigr).
	\end{align}
	
	From $\abs{\partial_p V}^2 \le C_V V$ we deduce that $\abs{\partial_p V^a} \le a \sqrt{C_V} V^{a-\frac{1}{2}}$ so that the application $\bigl(p \mapsto V^a(q,p)\bigr)$ is Lipschitz for every $a \le \frac{1}{2}$. 	
	Hence, if $a \le \frac{1}{2}$ we have \eqref{ac-Va-amontrer} with $A_{n-1} = \sqrt{V}(X_{n-1})$ and $B_n = 1 + \abs{U_n}$. If $a > \frac{1}{2}$, by Taylor's formula we have 
	\begin{align}
		V^a(X_n) - V^a(Q_n, P_{n-1}) & = \psca{\partial_p V^a(Q_n, \bar{p}_n), \Delta P_n}, \notag \\
		& \le a \sqrt{C_V} \bigl(\sqrt{V}(Q_n, \bar{p}_n)\bigr)^{2a-1} \abs{\Delta P_n}, \label{ineg-01}
	\end{align}
	with $\bar{p}_n \in (P_{n-1}, P_n)$ and since $2a-1 \le 1$  
	\begin{align}
		\bigl( \sqrt{V}(Q_n, \bar{p}_n) \bigr)^{2a-1} & \le \bigl( \sqrt{V}(Q_n, P_{n-1}) + C \abs{\Delta P_n} \bigr)^{2a - 1}, \notag \\
		& \le \bigl( \sqrt{V}(Q_n, P_{n-1}) \bigr)^{2a - 1} + C \abs{\Delta P_n}^{2a - 1}. \label{ineg-02}
	\end{align}
	Plugging \eqref{ineg-02} in \eqref{ineg-01} we get 
	\begin{multline*} 
		\abs{V^a(X_n) - V^a(Q_n, P_{n-1})} \\ \le C \sqrt{\gamma_n} \bigl(\sqrt{V}(Q_n, P_{n-1}) \bigr)^{2a-1} \sqrt{V}(X_{n-1}) \bigl(1 + \abs{U_n} \bigr) \\
		+ C \gamma_n^a \bigl(\sqrt{V}(X_{n-1}) \bigr)^{2a} \bigr(1+\abs{U_n} \bigr)^{2a}
	\end{multline*}
	Since $V \ge 1$ and $2a \ge 1$, it is easy to check that \eqref{ac-Va-amontrer} is satisfied with $A_{n-1} = \bigl( \sqrt{V}(Q_n, P_{n-1}) \bigr)^{2a-1} \bigl( \sqrt{V}(X_{n-1}) \bigr)^{2a}$ and $B_n = 1 + \abs{U_n}^{2a}$.

	By \eqref{ac-Va-amontrer} and Taylor's formula applied to the convex function $\bigl(x \mapsto \exp(\lambda x) \bigr)$ between $V^a(X_n)$ and $V^a(Q_n, P_{n-1})$ we have 
	\begin{equation*}
		\abs{W(X_n) - W(Q_n, P_{n-1}))} \le \lambda C \sqrt{\gamma_n} \bigl(W(X_n) + W(Q_n, P_{n-1}) \bigr) A_{n-1} B_n, 
	\end{equation*}
	and then 
	\begin{multline} \label{hamilton-tension-preresult}
		\E \bigl[ \abs{W(X_n) - W(Q_n, P_{n-1})}^s \,|\, \F_{n-1} \bigr] \le \lambda^s C^s \gamma_n^{s/2} A^s_{n-1} \E \bigl[W^s(X_n) B_n^s \,|\, \F_{n-1} \bigr] \\ + \lambda^s C^s \gamma_n^{s/2} A^s_{n-1} W^s(Q_n, P_{n-1}) \E [B_1^s].
	\end{multline}
	As in the proof of Lemma~\ref{lem-tension-essquad} we prove that for every $\varepsilon > 0$ there exists $C_\varepsilon > 0$ such that
	\begin{equation} \label{hamilton-tension-fin}
		\E \bigl[ W^s(X_n) B_n^s \,|\, \F_{n-1} \bigr] \le \frac{1}{1+\varepsilon} \exp \bigl(\lambda s(1+\varepsilon) V^a(X_{n-1}) \bigr) + C_\varepsilon.
	\end{equation}
	To complete the proof we apply Lemma~\ref{lem-hamilton-essentiel} which gives 
	\begin{equation*}
		V(Q_n, P_{n-1}) \le \bigl(1 + \sqrt{2 \delta C_V} \sqrt{\tilde{\gamma}_n} \bigr) V(X_{n-1}) + \tilde{\gamma}_n \beta,
	\end{equation*}
	and the concavity of $(x \mapsto x^a)$ implies 
	\begin{equation*}
		V^a(Q_n, P_{n-1}) \le \bigl(1 + a \sqrt{2 \delta C_V} \sqrt{\tilde{\gamma}_n} \bigr) V^a(X_{n-1}) + \tilde{\gamma}_n a \beta.
	\end{equation*}
	There exists $n_0 > 0$ such that for every $n \ge n_0$, $\lambda s (1 + a \sqrt{2\delta C_V} \sqrt{\gamma_n}) \le \bar{\lambda} < \lambda_0$, and thus there exists $K > 0$ such that 
	\begin{equation*}
		W^s(Q_n, P_{n-1}) \le K \exp \bigl(\bar{\lambda} V^a(X_{n-1}) \bigr).
	\end{equation*}
	Plugging this and \eqref{hamilton-tension-fin} in \eqref{hamilton-tension-preresult} and using the condition \eqref{hypothese-theta2} and Lemma~\ref{lem-clef} we obtain the result.
\qed \end{pf*}
\section{Numerical experiments}
The aim of this section is to validate our scheme numerically. 
We consider two problems: the Lorenz equations under external random excitation and a perturbed Hamiltonian system, which illustrate results in Sections~4 and 5 respectively. In both cases, we compute $\nu^\eta_n(f)$ for a given function $f$ using our scheme and we compare it with $\E \bigl[f(\bar{X}^T_h) \bigr]$ where $(\bar{X}^T_h)$ is a discretization scheme with constant step $h$ (an Euler scheme or an implicit Euler scheme). The approximation of $\E \bigl[f(\bar{X}^T_h) \bigr]$ is given by a Monte-Carlo procedure with 10000 paths. 
We give a representation of the stochastic sequence $(\tilde{\gamma}_n)_{n \ge 0}$ in Fig. \ref{fig-gamma-tilde-lorenz} and \ref{fig-gamma-tilde-3dof}. 

The programs are in C using BLAS/LAPACK (see {\tt http://www.netlib.org}) for linear algebra routine and the GSL library (see {\tt http://www.gnu.org}). In particular, the approximation of the fixed point needed in the implicit Euler scheme is done by the function {\tt gsl\_multiroot\_solver}. The random generator is a Mersenne twister generator of period $2^{19937}-1$ taken from {\tt http://www.math.sci.hiroshima-u.ac.jp/\~{}m-mat/MT/emt.html}.

All simulations are achieved using a Gaussian white noise for $(U_n)_{n \ge 1}$. The deterministic part $\gamma_n$ of the implemented step sequence is $\gamma_n = \gamma_0 n^{-1/3}$ and the weight sequence used is $(\eta_n)_{n \ge 1} = (\gamma_n)_{n \ge 1}$. These choices are motivated by the study of the rate of convergence (see \cite{these}).

\subsection{Lorenz equations}
It is a dissipative problem which is related to Section~4. The equations are 
\begin{equation*}
	\left\{ \begin{aligned}
	    \d x_t &= 10 (y_t - x_t) \d t + \d W^1_t, \\
	    \d y_t &= (28 x_t - y_t - x_t z_t) \d t + \d W^2_t, \\
	    \d z_t &= (x_t y_t -  (8/3) z_t) \d t
	\end{aligned} \right.
\end{equation*}
and a Lyapounov function for this system is $V(u)=\abs{u}^2 + 1$. The function $\chi$ used for simulations is defined for every $u \in \Rset^3$ by $\chi(u) = \frac{2 V(u)}{\abs{b(u)}^2 \vee 1}$. 

\begin{figure}[!ht]
\begin{center}
{\small %
\begin{picture}(0,0)%
\includegraphics{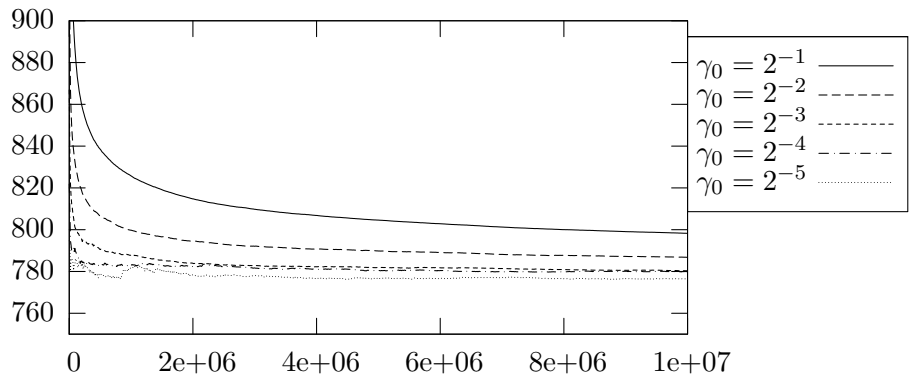}%
\end{picture}%
\begingroup
\setlength{\unitlength}{0.0200bp}%
\begin{picture}(18000,7559)(0,0)%
\put(1650,1494){\makebox(0,0)[r]{\strut{} 760}}%
\put(1650,2282){\makebox(0,0)[r]{\strut{} 780}}%
\put(1650,3070){\makebox(0,0)[r]{\strut{} 800}}%
\put(1650,3858){\makebox(0,0)[r]{\strut{} 820}}%
\put(1650,4646){\makebox(0,0)[r]{\strut{} 840}}%
\put(1650,5434){\makebox(0,0)[r]{\strut{} 860}}%
\put(1650,6222){\makebox(0,0)[r]{\strut{} 880}}%
\put(1650,7010){\makebox(0,0)[r]{\strut{} 900}}%
\put(1925,550){\makebox(0,0){\strut{} 0}}%
\put(4255,550){\makebox(0,0){\strut{} 2e+06}}%
\put(6585,550){\makebox(0,0){\strut{} 4e+06}}%
\put(8915,550){\makebox(0,0){\strut{} 6e+06}}%
\put(11245,550){\makebox(0,0){\strut{} 8e+06}}%
\put(13575,550){\makebox(0,0){\strut{} 1e+07}}%
\put(13575,6160){\makebox(0,0)[l]{\strut{} $\gamma_0=2^{-1}$}}%
\put(13575,5610){\makebox(0,0)[l]{\strut{} $\gamma_0=2^{-2}$}}%
\put(13575,5060){\makebox(0,0)[l]{\strut{} $\gamma_0=2^{-3}$}}%
\put(13575,4510){\makebox(0,0)[l]{\strut{} $\gamma_0=2^{-4}$}}%
\put(13575,3960){\makebox(0,0)[l]{\strut{} $\gamma_0=2^{-5}$}}%
\end{picture}%
\endgroup
}\caption{\label{nu_lorenz} One path of $\bigl(\nu^\eta_n(f)\bigr)_{1 \le n \le 10^7}$ for different value of $\gamma_0$.}
\end{center}
\end{figure}

The stochastic step sequence used for our scheme is thus 
\begin{equation}
	\tilde{\gamma}_0 = \gamma_0, \quad \text{and} \quad \forall n \ge 1, \quad \tilde{\gamma}_n = \bigl( \gamma_0 n^{-\frac{1}{3}} \bigr) \wedge \Bigl(\frac{2 V(X_{n-1})}{\abs{b(X_{n-1})}^2 \vee 1} \Bigr).
\end{equation}
The weight sequence $(\eta_n)_{n \ge 1}$ is equal to $(\gamma_n)_{n \ge 1}$ and we compute for $f(x) = \abs{x}^2$ and $n \le 10^7$ 
\begin{equation}
	\nu^\eta_n(f) = \frac{1}{\sum_{k=1}^n k^{-\frac{1}{3}}} \sum_{k=1}^n k^{-\frac{1}{3}} f(X_{k-1}).
\end{equation}
The result is given by the figure \ref{nu_lorenz}.
The scheme seems to have a better behaviour when $\gamma_0$ is equals to $2^{-4}$ or $2^{-5}$. For bigger values of $\gamma_0$, the rate of convergence is poor but the Lorenz problem is a difficult numerical problem and the parameter $\gamma_0$ is hard to fix. Other numerical methods have the same problem. The important point to note here is that the scheme does not explode (for any $\gamma_0$) and appears convergent to the same limit. 

\begin{figure}[!ht]
\begin{center}
{\small %
\begin{picture}(0,0)%
\includegraphics{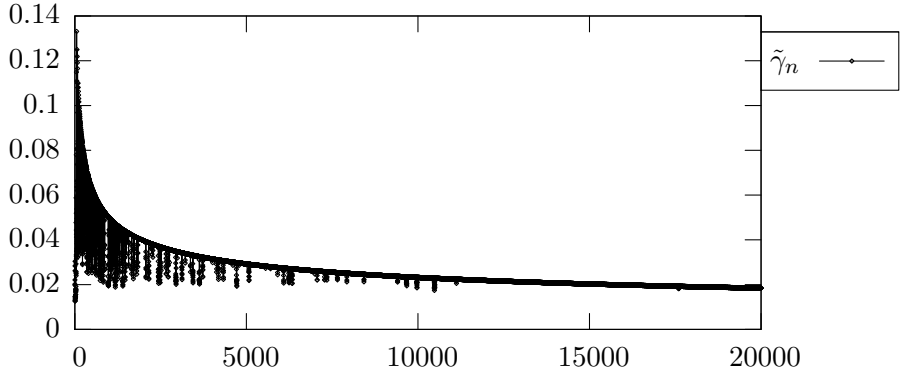}%
\end{picture}%
\begingroup
\setlength{\unitlength}{0.0200bp}%
\begin{picture}(18000,7559)(0,0)%
\put(1925,1100){\makebox(0,0)[r]{\strut{} 0}}%
\put(1925,1944){\makebox(0,0)[r]{\strut{} 0.02}}%
\put(1925,2789){\makebox(0,0)[r]{\strut{} 0.04}}%
\put(1925,3633){\makebox(0,0)[r]{\strut{} 0.06}}%
\put(1925,4477){\makebox(0,0)[r]{\strut{} 0.08}}%
\put(1925,5321){\makebox(0,0)[r]{\strut{} 0.1}}%
\put(1925,6166){\makebox(0,0)[r]{\strut{} 0.12}}%
\put(1925,7010){\makebox(0,0)[r]{\strut{} 0.14}}%
\put(2200,550){\makebox(0,0){\strut{} 0}}%
\put(5432,550){\makebox(0,0){\strut{} 5000}}%
\put(8665,550){\makebox(0,0){\strut{} 10000}}%
\put(11897,550){\makebox(0,0){\strut{} 15000}}%
\put(15130,550){\makebox(0,0){\strut{} 20000}}%
\put(15130,6160){\makebox(0,0)[l]{\strut{} $\tilde{\gamma}_n$}}%
\end{picture}%
\endgroup
}\caption{\label{fig-gamma-tilde-lorenz} Stochastic step sequence $(\tilde{\gamma}_n)_{n \ge 1}$.}
\end{center}
\end{figure}

Figure \ref{fig-gamma-tilde-lorenz} gives a representation of the stochastic step sequence $(\tilde{\gamma}_n)_{n \ge 1}$ when $\gamma_0 = 0.5$. We show that the bigger $n$ is and the less the stochastic part $\chi(X_{n-1})$ is used. For this path, after $20\,000$ iterations we have $\tilde{\gamma}_n = \gamma_n$ (at least until $n = 10^7$). Moreover before the $20\,000$th iteration the event $\{\chi(X_{n-1}) < \gamma_n\}$ occurs only 924 times (4.62\% of time).
 
\begin{figure}[!ht]
\begin{center}
	\subfigure[\label{nu_h_lorenz_euler} Euler scheme: $\bar{X}^T_{k+1,h} = \bar{X}^T_{k,h} + h b(\bar{X}^T_{k,h}) + \sqrt{h} \sigma U_{k+1}$]{\small %
\begin{picture}(0,0)%
\includegraphics{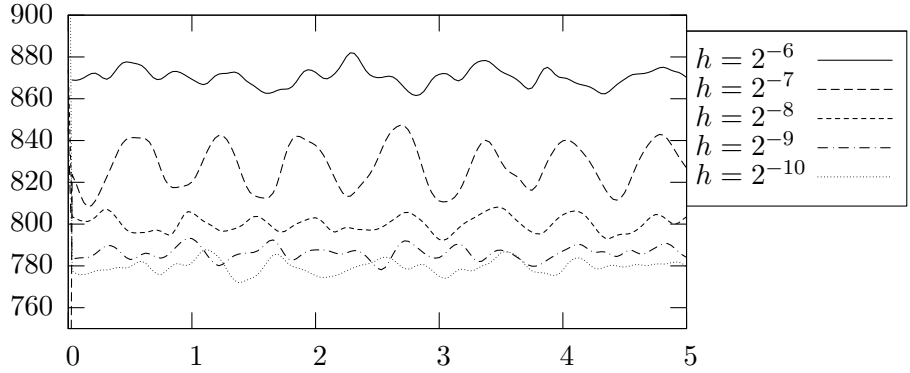}%
\end{picture}%
\begingroup
\setlength{\unitlength}{0.0200bp}%
\begin{picture}(18000,7559)(0,0)%
\put(1650,1494){\makebox(0,0)[r]{\strut{} 760}}%
\put(1650,2282){\makebox(0,0)[r]{\strut{} 780}}%
\put(1650,3070){\makebox(0,0)[r]{\strut{} 800}}%
\put(1650,3858){\makebox(0,0)[r]{\strut{} 820}}%
\put(1650,4646){\makebox(0,0)[r]{\strut{} 840}}%
\put(1650,5434){\makebox(0,0)[r]{\strut{} 860}}%
\put(1650,6222){\makebox(0,0)[r]{\strut{} 880}}%
\put(1650,7010){\makebox(0,0)[r]{\strut{} 900}}%
\put(1925,550){\makebox(0,0){\strut{} 0}}%
\put(4255,550){\makebox(0,0){\strut{} 1}}%
\put(6585,550){\makebox(0,0){\strut{} 2}}%
\put(8915,550){\makebox(0,0){\strut{} 3}}%
\put(11245,550){\makebox(0,0){\strut{} 4}}%
\put(13575,550){\makebox(0,0){\strut{} 5}}%
\put(13575,6160){\makebox(0,0)[l]{\strut{} $h=2^{-6}$}}%
\put(13575,5610){\makebox(0,0)[l]{\strut{} $h=2^{-7}$}}%
\put(13575,5060){\makebox(0,0)[l]{\strut{} $h=2^{-8}$}}%
\put(13575,4510){\makebox(0,0)[l]{\strut{} $h=2^{-9}$}}%
\put(13575,3960){\makebox(0,0)[l]{\strut{} $h=2^{-10}$}}%
\end{picture}%
\endgroup
} \\	\subfigure[\label{nu_h_lorenz_back} Implicit Euler scheme: $\bar{X}^T_{k+1,h} = \bar{X}^T_{k,h} + h b(\bar{X}^T_{k+1,h}) + \sqrt{h} \sigma U_{k+1}$]{\small %
\begin{picture}(0,0)%
\includegraphics{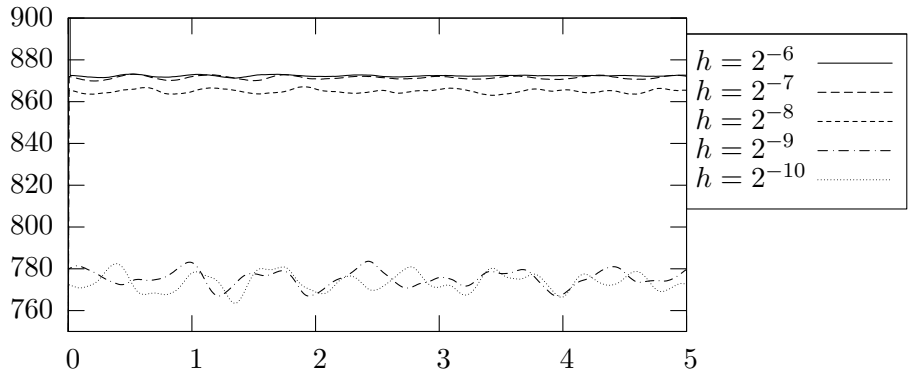}%
\end{picture}%
\begingroup
\setlength{\unitlength}{0.0200bp}%
\begin{picture}(18000,7559)(0,0)%
\put(1650,1494){\makebox(0,0)[r]{\strut{} 760}}%
\put(1650,2282){\makebox(0,0)[r]{\strut{} 780}}%
\put(1650,3070){\makebox(0,0)[r]{\strut{} 800}}%
\put(1650,3858){\makebox(0,0)[r]{\strut{} 820}}%
\put(1650,4646){\makebox(0,0)[r]{\strut{} 840}}%
\put(1650,5434){\makebox(0,0)[r]{\strut{} 860}}%
\put(1650,6222){\makebox(0,0)[r]{\strut{} 880}}%
\put(1650,7010){\makebox(0,0)[r]{\strut{} 900}}%
\put(1925,550){\makebox(0,0){\strut{} 0}}%
\put(4255,550){\makebox(0,0){\strut{} 1}}%
\put(6585,550){\makebox(0,0){\strut{} 2}}%
\put(8915,550){\makebox(0,0){\strut{} 3}}%
\put(11245,550){\makebox(0,0){\strut{} 4}}%
\put(13575,550){\makebox(0,0){\strut{} 5}}%
\put(13575,6160){\makebox(0,0)[l]{\strut{} $h=2^{-6}$}}%
\put(13575,5610){\makebox(0,0)[l]{\strut{} $h=2^{-7}$}}%
\put(13575,5060){\makebox(0,0)[l]{\strut{} $h=2^{-8}$}}%
\put(13575,4510){\makebox(0,0)[l]{\strut{} $h=2^{-9}$}}%
\put(13575,3960){\makebox(0,0)[l]{\strut{} $h=2^{-10}$}}%
\end{picture}%
\endgroup
}	\caption{\label{nu_h_lorenz} Computation of $\E \bigl[f(\bar{X}^T_h)\bigr]$ by a Monte-Carlo procedure on 10000 paths and for $T = 5$}
\end{center}
\end{figure}

We compare our results with the approximation of $\E \bigl[f(\bar{X}^T_h)\bigr]$  where $(\bar{X}^T_h)$ is a regular Euler scheme (figure \ref{nu_h_lorenz_euler}) or a implicit Euler scheme (figure \ref{nu_h_lorenz_back}). We represent the results only for $h \le 2^{-6}$ because for bigger values the empirical expectation (based on $10\,000$ paths) explodes whith the regular Euler scheme. For the implicit scheme the expectation remains bounded but the behaviour is very poor when $h \le 2^{-8}$. The parameters $h$, $T$ and the number of paths used for the Monte-Carlo procedure are hard to fix. 

For the implicit Euler scheme (figure \ref{nu_h_lorenz_back}), the jump between $h = 2^{-8}$ and $h = 2^{-9}$ is due to the fact that the scheme remains trapped in the neighborhood of only one attractor (Lorenz equation has two attractors). This behavior does not occur with the regular Euler scheme and with our scheme. 

\subsection{Perturbed Hamiltonian system}
The second example is a perturbed Hamiltonian system derived from a multidimensional linear oscillator under external random excitation. It is a 3-DOF (degree of freedom) system studied by Ibrahim and Li (see \cite{ibrahim-li}) and Soize (see \cite{soize} chap. XIII.6). We have thus the following equation in $\Rset^6$
\begin{equation*}
	\left\{ \begin{aligned}
	    \d q_t &= \partial_p H(q_t, p_t) \d t, \\
	    \d y_t &= - \partial_q H(q_t, p_t) \d t - f_0 D_0 \partial_p H(q_t, p_t) \d t + g_0 S_0 \d W_t
	\end{aligned} \right.
\end{equation*}
where $H(q_t, p_t) = \frac{1}{2} \psca{M(q)^{-1} p, p} + \frac{1}{2} \psca{K_0 q, q}$ and 
\begin{equation*}
	D_0 = S_0 = \begin{pmatrix} 
	    1 & 0 & 0 \\
	    0 & 0 & 0 \\
	    0 & 0 & 0 
	\end{pmatrix}, \quad 
	K_0 = \begin{pmatrix}
	    1.3 & 0 & 0 \\
	    0 & 0.154 & 0 \\
	    0 & 0 & 0.196 
	\end{pmatrix},  
	\quad g_0 = 0.5, \quad f_0 = 0.9965,
\end{equation*}
and 
\begin{equation*}
	M(q) = \begin{pmatrix}
	    1.3 & v_2 + 0.3 & v_3  \\
	    v_2 + 0.3 & v_2^2 + 2 v_2 + 0.314 & v_2 v_3 + v_3 + 0.0375 \\
	    v_3 & v_2 v_3 + v_3 + 0.0375 & v_3^2 + 0.1
	\end{pmatrix}, \quad 
\end{equation*}
with $v_2(q) = - 1.61 (0.375 q_2 + q_3)$ and $v_3(q) = - 1.61 q_2$. This numerical example is taken from \cite{soize} (page 257--264). This is the first damping model case with the external excitation applied to DOF 1 and the system parameter equal to 0.7. In this case we have an analytic expression of the density of the invariant measure and we can calculate the mean-square response for the DOF 1. We consider the function $f_1:(q,p) \mapsto q_1^2$ and we have 
\begin{equation*}
	\int_{\Rset^6} f_1(x) \nu(\d x) \simeq 0.0965.
\end{equation*}

\begin{figure}[!ht]
\begin{center} 
{\small %
\begin{picture}(0,0)%
\includegraphics{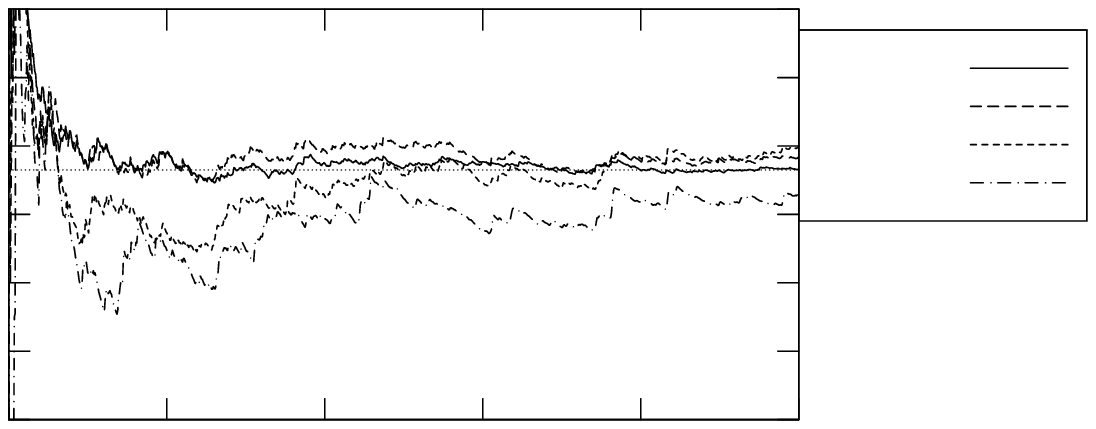}%
\end{picture}%
\begingroup
\setlength{\unitlength}{0.0200bp}%
\begin{picture}(18000,7559)(0,0)%
\put(1925,1100){\makebox(0,0)[r]{\strut{} 0.06}}%
\put(1925,2085){\makebox(0,0)[r]{\strut{} 0.07}}%
\put(1925,3070){\makebox(0,0)[r]{\strut{} 0.08}}%
\put(1925,4055){\makebox(0,0)[r]{\strut{} 0.09}}%
\put(1925,5040){\makebox(0,0)[r]{\strut{} 0.1}}%
\put(1925,6025){\makebox(0,0)[r]{\strut{} 0.11}}%
\put(1925,7010){\makebox(0,0)[r]{\strut{} 0.12}}%
\put(2200,550){\makebox(0,0){\strut{} 0}}%
\put(4475,550){\makebox(0,0){\strut{} 200000}}%
\put(6750,550){\makebox(0,0){\strut{} 400000}}%
\put(9025,550){\makebox(0,0){\strut{} 600000}}%
\put(11300,550){\makebox(0,0){\strut{} 800000}}%
\put(13575,550){\makebox(0,0){\strut{} 1e+06}}%
\put(13575,6160){\makebox(0,0)[l]{\strut{} $\gamma_0 = 2^{-1}$}}%
\put(13575,5610){\makebox(0,0)[l]{\strut{} $\gamma_0 = 2^{-2}$}}%
\put(13575,5060){\makebox(0,0)[l]{\strut{} $\gamma_0 = 2^{-3}$}}%
\put(13575,4510){\makebox(0,0)[l]{\strut{} $\gamma_0 = 2^{-4}$}}%
\end{picture}%
\endgroup
}\caption{\label{nu_3dof} One path of $\bigl(\nu^\eta_n(f)\bigr)_{1 \le n \le 10^6}$ for different value of $\gamma_0$.}
\end{center}
\end{figure}

\begin{figure}[!ht]
\begin{center}
{\small %
\begin{picture}(0,0)%
\includegraphics{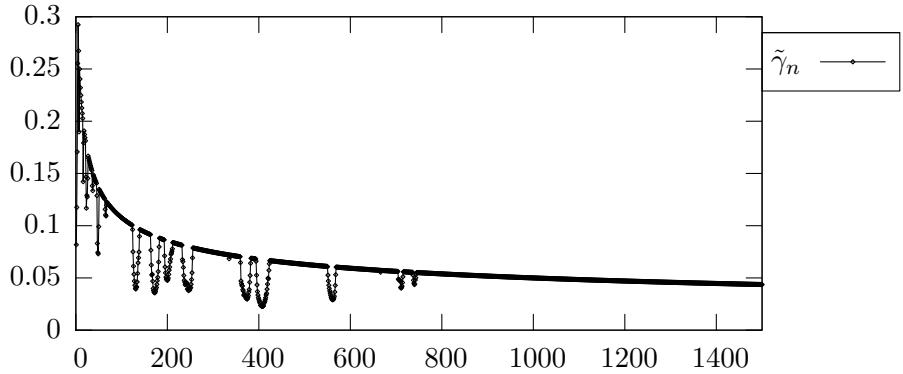}%
\end{picture}%
\begingroup
\setlength{\unitlength}{0.0200bp}%
\begin{picture}(18000,7559)(0,0)%
\put(1925,1100){\makebox(0,0)[r]{\strut{} 0}}%
\put(1925,2085){\makebox(0,0)[r]{\strut{} 0.05}}%
\put(1925,3070){\makebox(0,0)[r]{\strut{} 0.1}}%
\put(1925,4055){\makebox(0,0)[r]{\strut{} 0.15}}%
\put(1925,5040){\makebox(0,0)[r]{\strut{} 0.2}}%
\put(1925,6025){\makebox(0,0)[r]{\strut{} 0.25}}%
\put(1925,7010){\makebox(0,0)[r]{\strut{} 0.3}}%
\put(2200,550){\makebox(0,0){\strut{} 0}}%
\put(3924,550){\makebox(0,0){\strut{} 200}}%
\put(5648,550){\makebox(0,0){\strut{} 400}}%
\put(7372,550){\makebox(0,0){\strut{} 600}}%
\put(9096,550){\makebox(0,0){\strut{} 800}}%
\put(10820,550){\makebox(0,0){\strut{} 1000}}%
\put(12544,550){\makebox(0,0){\strut{} 1200}}%
\put(14268,550){\makebox(0,0){\strut{} 1400}}%
\put(15130,6160){\makebox(0,0)[l]{\strut{} $\tilde{\gamma}_n$}}%
\end{picture}%
\endgroup
}\caption{\label{fig-gamma-tilde-3dof} One representation of the stochastic step sequence $(\tilde{\gamma}_n)_{n \ge 1}$.}
\end{center}
\end{figure}

The stochastic sequence used in the following simulations is defined by $\tilde{\gamma}_0 = \gamma_0$ and $\forall n \ge 1$ 
\begin{equation*}
	\tilde{\gamma}_n = \bigl( \gamma_0 n^{-1/3} \bigr) \wedge \Bigl( \frac{1}{\abs{b(X_{n-1})}^2 \vee 1} \Bigr).
\end{equation*}
The results for different value of $\gamma_0$ are given in figure \ref{nu_3dof}. The convergence seems better when $\gamma_0$ is big. Our scheme behaves very well and a representation of the stochastic step sequence $(\tilde{\gamma}_n)_{n \ge 1}$ is given in figure \ref{fig-gamma-tilde-3dof}. 

We do not represent the approximation of $\E \bigl[f(\bar{X}^T_h)\bigr]$ where $(\bar{X}^T_h)$ is the explicit Euler scheme because the empirical expectation (based on $10\,000$ paths) explodes for $h < 2^{-3}$.

\begin{figure}[!ht]
\begin{center}
{\small %
\begin{picture}(0,0)%
\includegraphics{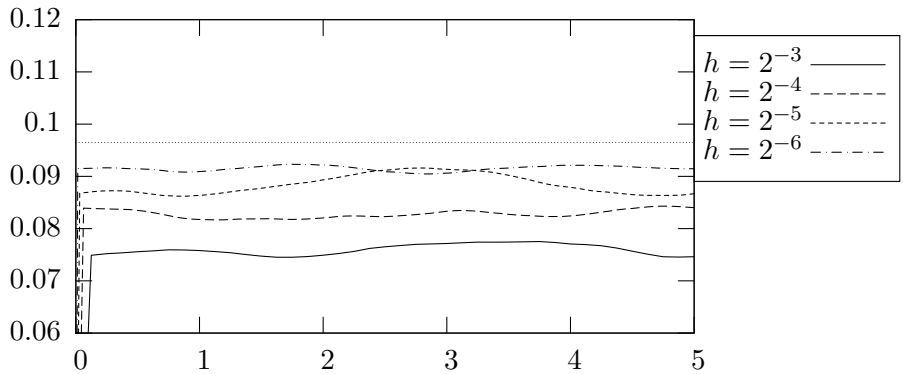}%
\end{picture}%
\begingroup
\setlength{\unitlength}{0.0200bp}%
\begin{picture}(18000,7559)(0,0)%
\put(1925,1100){\makebox(0,0)[r]{\strut{} 0.06}}%
\put(1925,2085){\makebox(0,0)[r]{\strut{} 0.07}}%
\put(1925,3070){\makebox(0,0)[r]{\strut{} 0.08}}%
\put(1925,4055){\makebox(0,0)[r]{\strut{} 0.09}}%
\put(1925,5040){\makebox(0,0)[r]{\strut{} 0.1}}%
\put(1925,6025){\makebox(0,0)[r]{\strut{} 0.11}}%
\put(1925,7010){\makebox(0,0)[r]{\strut{} 0.12}}%
\put(2200,550){\makebox(0,0){\strut{} 0}}%
\put(4530,550){\makebox(0,0){\strut{} 1}}%
\put(6860,550){\makebox(0,0){\strut{} 2}}%
\put(9190,550){\makebox(0,0){\strut{} 3}}%
\put(11520,550){\makebox(0,0){\strut{} 4}}%
\put(13850,550){\makebox(0,0){\strut{} 5}}%
\put(13850,6160){\makebox(0,0)[l]{\strut{} $h=2^{-3}$}}%
\put(13850,5610){\makebox(0,0)[l]{\strut{} $h=2^{-4}$}}%
\put(13850,5060){\makebox(0,0)[l]{\strut{} $h=2^{-5}$}}%
\put(13850,4510){\makebox(0,0)[l]{\strut{} $h=2^{-6}$}}%
\end{picture}%
\endgroup
}\caption{\label{nu_h_3dof} Computation of $\E \bigl[f(\bar{X}^T_h)\bigr]$ by a Monte-Carlo procedure on $10\,000$ paths and for $T = 5$ where $\bar{X}^T_h$ is a implicit Euler scheme.}
\end{center}
\end{figure}

Figure~\ref{nu_h_3dof} gives results using the approximation of $\E \bigl[f(\bar{X}^T_h)\bigr]$ where $(\bar{X}^T_h)$ is the implicit Euler scheme. 

\bibliography{biblio}
\bibliographystyle{plain}

\end{document}